\documentclass{amsart}

\usepackage{amssymb,latexsym,amsfonts,amsmath}
\usepackage{graphicx,psfrag,epsfig,epsf}
\usepackage{eso-pic}
\usepackage{graphicx}
\usepackage{color}
\usepackage{tikz}
\usetikzlibrary{automata}
\usepackage[vlined,ruled]{algorithm2e}

\usepackage{amsfonts}
\usepackage{amssymb,latexsym,amsfonts,amsmath}
\usepackage{graphicx}

\newtheorem{theorem}{Theorem}[section]
\newtheorem{lemma}[theorem]{Lemma}

\newtheorem{proposition}[theorem]{Proposition}
\newtheorem{corollary}[theorem]{Corollary}

\theoremstyle{definition}
\newtheorem{definition}[theorem]{Definition}

\theoremstyle{remark}
\newtheorem{remark}[theorem]{Remark}
\numberwithin{equation}{section}

\SetKwInput{KwAssumes}{Assumes}

\def\re{{\rm I\!R}}
\def\na{{\rm I\!N}}
\def\co{{\rm I\!\!\!C}}

\def\cal{\mathcal}

\def\be{\begin{equation}}
\def\ee{\end{equation}}
\def\bea{\begin{eqnarray}}
\def\eea{\end{eqnarray}}

\def\bmat{\left[\begin{array}}
\def\emat{\end{array}\right]}

\begin{document}

\begin{abstract}                          % Abstract of not more than 200 words.
We show a general method allowing the solution calculation, in the form of a power series,  
for a very large class of  nonlinear  
Ordinary Differential Equations (ODEs), namely the real analytic $\sigma\pi$-ODEs (and, more in general, the real analytic 
$\sigma\pi$-{\it reducible} ODEs) in many indeterminates, 
characterized by an ODE-function given by generalized polynomials of the indeterminates and their derivatives, i.e. functions formally polynomial with exponents, though the exponent can be any real number.  
The solution method consists in reducing 
the ODE to a certain canonical quadratic ODE, with a larger number of indeterminates, whose solutions include, as sub-solutions, 
the original solutions, and for which a recursive formula is shown, giving all coefficients of the solution power series directly from the ODE parameters.  The reduction method  is named {\it exact quadratization} and was formerly introduced in another our article, where 
we considered explicit ODEs only.  In the present paper, which is self-contained to a large extent, 
we review and complete the theory of exact quadratization by solving issues, such as for instance the {\it globality} of the quadratization, that had remained open, and also extend it to the more general case of an implicitly defined ODE. Finally, we argue that the result can be seen as 
a partial solution of a differential version of the 22nd Hilbert's problem. 
\end{abstract}

\title{On the Solution Calculation of  Nonlinear Ordinary Differential Equations via Exact Quadratization}
%\thanks{This material is based upon work supported in part by ARO
%MURI Grant number W911NF0910553, and in part by the Center of Excellence for Research DEWS, University of L'Aquila, Italy}

\author[Francesco Carravetta]
{Francesco Carravetta$^{1}$
%, Shaunak D. Bopardikar$^{2}$, Joao P. Hespanha$^{2}$, Maria D. Di Benedetto$^{1}$
}
\address{$^{1}$
Istituto di Analisi dei
  Sistemi ed Informatica ``Antonio Ruberti'', Consiglio Nazionale della Ricerche, Via dei Taurini 19, Roma,
  Italy.
}
\email{francesco.carravetta@iasi.cnr.it}
%\address{$^{2}$
%Center for Control, Dynamical Systems and Computation\\
 % University of California at Santa Barbara, CA 93106, USA}
%\email{ \{sdbopardikar,hespanha\}@ece.ucsb.edu}

\maketitle

\section{Introduction}

\label{sc:int}

Let us consider the following quadratic ordinary differential equation (ODE) in $\re^m$, written component-wise ($i=1,\ldots,m$): 
\be
\dot x_i=\sum_{j=1}^m v_{i,j}x_ix_j,\label{Quad1}
\ee
where the 'time' $t\in\cal T\subset \re$, $\cal T$ an open interval, (resp: $x\in\re^m$) is the independent (resp: dependent) variable, the dot 
denotes time-derivative, $x_i$ denotes the $i$th component of the vector $x$, and $v_{i,j}$ are real functions of $t$ analytic on $\cal T$. 
Whatever we fix a point $(x,s)\in \re^n\times \cal T$ we know that there exists an unique solution $t\mapsto x(t)$ of (\ref{Quad1}) defined and analytic on some interval ${\bf I}_{x,s}\subset \cal T$, which in general depends of the point $(x,s)$, with $s\in {\bf I}_{x,s}$ and $x(s)=x$. Moreover, since the solution $x(t)$ is analytic, for any $t_0\in {\bf I}_{x,s}$  it can be written, component-wise, as the power series:  
\be
x_i(t)=\sum_{k=0}^{\infty}c_k(i) {(t-t_0)^k\over k!},\qquad c_k(i)=x_i^{(k)}(t_0), \label{SPSDint}
\ee
where $x^{(k)}$ denotes the $k$th time-derivative, and (\ref{SPSDint}) holds at least whereby all series, for $i=1,\ldots,n$, converge, which surely occurs for all  $t\in {\bf I}(t_0,r) \subset 
{\bf I}_{x,s}$, an open interval centered in $t_0$ with ray $r\in \re^+\cup \{+\infty\}$. \\

Now, let us pose the following question: \\

{\it (q)} is there a recursive (in $k$) formula giving the coefficients $c_k(i)$ of the power series (\ref{SPSDint}) for all $k\in \na$ directly from the coefficients $v_{i,j}$?\\

The above question has an interest going far beyond the single case study, insofar one is able to show that the quadratic ODE (\ref{Quad1}) {\it is not} really a {\it particular} differential equation, and it is in fact a sort of {\it nonlinear paradigm}, a basic case to which a large class of analytic ODE can be always {\it reduced}, in the sense we are going to explain below. \\

We say that a first order ODE of the type: 
\be
\dot x=f(x,t),\label{S1}
\ee
where $f:\cal A\to \re^n$, is a function defined on the open set  $\cal A\subset \re^{n+1}$ and locally Lipschitz therein, 
is {\it globally, exactly} (adjectives that hereinafter we omit  whether clear from context) {\it quadratizable}, or that {\it undergoes a (global, exact) quadratization}, if there exist an integer $m\ge n$, and locally Lipschitz functions 
$v_{i,j}: \cal T\to \re$ -- with $\cal T$ an open interval of $\re$ and $i,j=1,\ldots,m$ --  such that 
every solution is a {\it local sub-solution}\footnote{We state in a moment in which sense we talk about a {\it sub-solution}.} of a solution of a  quadratic ODE of the type (\ref{Quad1}), which is then accordingly referred to as {\it a  quadratization} of (\ref{S1}).

Before going more in depth, we clarify some points of the definition given above.  
We define the {\it maximal domain}, $\cal A\subset  \re^{n+1}$, of (\ref{S1}) to be the largest (open) set such that for any $(x,t)\in\cal A$  
there exist a unique solution of (\ref{S1}) whose graph passes through $(x,t)$. We always understand that $\cal A$ is non void. 
If $(x,t)\in \cal A$, we denote by ${\bf I}_{x,t}\subset \re$ the largest (open) interval where 
the solution $s\mapsto x(s)$, such that $x(t)=x$ is defined, and identify the solution with the point $(x,t)$, or even with the sole $x$, when 
$t$ is understood, obvious, or un-essential \footnote{For any pair of points in $\cal A$ we can say they are equivalent if they are connected by the graph of a solution. Indeed, one readily verifies that this is an equivalence relation. This motivates the identification of the graph (and, hence, of the solution) with 
any point $(x,t)$ of the graph itself.}. 
We say that $(x,t)\in \cal A$ -- solution of the ODE (\ref{S1}) -- is a local (resp: global, or no adjective) sub-solution of  a second ODE, of the same type (\ref{S1}) but $m$-dimensional, with $m\ge n$,  if it is equal, point-wise in  
${\bf I}\subset {\bf I}_{x,t}$, with ${\bf I}$ an open sub-interval, (resp: point-wise in ${\bf I}_{x,t}$)  
to $n$ out of the $m$ components of a solution  of the second ODE. If we collect  the indices of the latter $n$ components in a set 
$\cal I^*\subset \{1,\ldots,m\}$ we can say that $(x,t)$ is a sub-solution of the second ODE {\it with respect to} $\cal I^*$.

Let us consider a map 
$\Phi:{\rm Dom}\{\Phi\}{\rm (open)}\subset {\cal A}\to \re^m$, satisfying the following property: for any $(x,t) \in {\rm Dom}\{\Phi\}$ there exists 
an open neighborhood $U\subset {\rm Dom}\{\Phi\}$ of $(x,t)$ such that, for any $(x',t')\in U$ belonging to the graph of the solution $(x,t)$ of (\ref{S1}),  $(\Phi(x',t'),t')$ belongs to the graph of the solution $(\Phi(x,t),t)$ of (\ref{Quad1}). We call such a map: {\it locally compatible } 
(with the ODEs (\ref{S1}) and (\ref{Quad1})), \footnote{In other words $\Psi$ preserves (locally) the equivalence relation between points of the same solution, 
i.e. maps (locally) solution of (\ref{S1}) into solution of (\ref{Quad1}). In a differential geometric framework (see for instance \cite{Bo75}) that  is  the well known property of {\it relatedness} of  the vector fields associated to the two ODEs.}  
and 
the ODE (\ref{S1}) will be said {\it (globally, exactly) 
quadratizable under $\Phi$} if every solution $(x,t)$ of (\ref{S1}) is a local sub-solution of the type 
$(\Phi(x^*,t^*),t^*)$ -- with $(x^*,t^*)\in {\rm Dom}\{\Phi\}$ -- of  a quadratic ODE of the type (\ref{Quad1}). 
If (\ref{S1}) is quadratizable under $\Phi$, we always understand in the following that $\Phi$ is a locally compatible map. 

If (\ref{S1}) is quadratizable under $\Phi$, and any of its solutions $(x,t)$ is an analytic function on ${\bf I}_{x,t}$ --- 
which is the case as $f$ in $(\ref{S1})$ is analytic in all of its arguments ---  
then any solution $(x,t)$ of (\ref{S1}) whose graph crosses ${\rm Dom}\{\Phi\}$,  say $(x(t^*),t^*)=(x^*,t^*)\in {\rm Dom}\{\Phi\}$,  
can be recovered (by analytic prolongation) by the 
solution $(\Phi(x^*,t^*),t^*)$ of (\ref{Quad1}).  In other words: we can {\it calculate} at least all solutions whose graph is not disjoint with the domain of $\Phi$, provided that we are able to calculate, even only locally, the solutions of a quadratic ODE of the type (\ref{Quad1}). 

In the paper \cite{Ca15} we have given a first answer to the question 
about which kind of ODEs can be exactly quadratized. 
The basic result (stated in a rather different setting than the one we use here)  has been that all ODEs  of the type (\ref{S1}) with $f$ given component-wise by: 
\bea
f_i(x,t)=\sum_{l=1}^{\nu_i} v_{i,l}(t) X_{i,l}(x); \quad X_{i,l}(x)=\prod_{j=1}^n x_j^{p^l_{i,j}},\label{SP-form1}
\eea
where $p^l_{i,j}$ are real exponents, $\nu_i\in\na$ (if $\nu_i=0$ then $f_i=0$) and $v_{i,l}$ are analytic functions all 
defined on some open real interval, 
namely $\cal T$, 
undergo a {\it dense} quadratization (we clarify in a moment what 'dense quadratization' means) under 
the following real analytic map\footnote{In \S II, it will be made clear that $\Phi$ is meromorphic on $\cal A$, but analytic on $\cal A\setminus \{x\in\re^n | x_i=0, i=1,\ldots,n\}={\rm Dom}\{\Phi\}$} $\Phi:{\rm Dom}\{\Phi\}\subset {\cal A}\to \re^m$ -- $\cal A\subset \re^n\times \cal T$ being the (open) analyticity domain of (\ref{SP-form1}),  and $m=\nu_1+\ldots+\nu_n$, thus $m\ge n$ -- which can be 
written component-wise as follows:
\be
\Phi_{i,l}(x_1,\ldots,x_n)=x_i^{-1} X_{i,l}. \label{ChV}
\ee  

In (\ref{SP-form1}) the function $f$  is a {\it formal polynomial} of $x$, that is to say: a polynomial writing where the exponents are allowed to be {\it any real number}. In \cite{Ca15} we named this kind of functions {\it $\sigma\pi$-functions}, and 
$\sigma\pi$-{\it ODE} the associated differential equation. 

A {\it dense} quadratization under a map $\Phi$ differs from the definition of (global) quadratization (under $\Phi$) given before in that it is not required that {\it every} solution-graph 
crosses ${\rm Dom}\{\Phi\}$, but rather that the set of all solution-graphs that do not cross it form a zero-measure set in $\re^{n+1}$, 
i.e. the set $\cal A\cap {\rm Dom}\{\Phi\}$ has zero-measure, or, which is the same, ${\rm Dom}\{\Phi\}$ is a {\it dense} subset in $\cal A$. Indeed, it can be at once inferred by (\ref{ChV}) that ${\rm Dom}\{\Phi\}$ is $\cal A\setminus S$ where 
$S$ denotes the union of the $n$ hyperplanes in $\re^{n+1}$ described by the equation $x_i=0$,  for $i=1,\ldots,n$, 
so,  what has been proven in \cite{Ca15} can be stated in other words as follows: 
every solution of a $\sigma\pi$-system, except possibly a sub-set whose curves all lye in some (among the first $n$)  
{\it coordinate hyperplane} of $\re^{n+1}$, is a local sub-solution of  an ODE of type (\ref{Quad1}) for some $m\ge n$ (with different coefficients $v_{i,l}$).

We can similarly say that (\ref{S1}) is $\sigma\pi$-{\it reducible} (resp: under a map $\Phi$), if any solution 
$(x,t)$ of it, is a sub-solution (resp: a sub-solution of the type $(\Phi(x,t),t)$) of the solution  of some  $\sigma\pi$-ODE (generally of a larger dimension). 
Obviously, the class of ODEs (\ref{S1}) that are exactly quadratizable includes all those $\sigma\pi$-reducible, and as shown in \cite{Ca15} there are indeed {\it many} of such functions $f$, for instance: {\it (i)} all finite compositions of 
$\sigma\pi$-functions, {\it (ii)} for $x\in \re$, the functions $e^x,\ln x$, and all the trigonometric functions $\sin x,\cos x, \tan x$, {\it (iii)} 
all finite compositions of functions from (i) and (ii). The problem of better characterize the {\it largest} set of functions that are densely quadratizable (or, equivalently: $\sigma\pi$-reducible) in terms of the more usual concepts of {\it algebraic} or {\it differentially algebraic} function, was therein only sketched, and we'll not be concerned here either with  this important issue, but save it for another occasion.  
It should be stressed that, 
on account that the class of $\sigma\pi$-ODEs  is by itself quite large, 
the class of all quadratizable ODEs, 
despite based on \cite{Ca15} has still vaguely defined outlines, 
appears to be 
{\it huge} as it is, because it includes all ODEs that one encounter in practice, typically defined through finite compositions of elementary, algebraic or transcendental functions of the type we have seen just above. With that being said,  the property for an ODE of being exactly quadratizable 
appears not ensuing given 
some (substantially) {\it restrictive} condition, and in fact to be a somewhat universal property of ordinary differential equations 
{\it of the type} (\ref{S1}), i.e. {\it explicitable} with respect to the first order derivative.

The aim of the present paper is, first, to give an answer to some issues that remained unsolved in \cite{Ca15}. The first of these issues 
is about whether 
an analytic $\sigma\pi$-ODE is {\it globally} quadratizable (and not just {\it densely} quadratizable). 
In other words, this amounts to check whether solutions $(x,t)\not\in {\rm Dom}\{\Phi\}$ are still local sub-solution of an ODE of type (\ref{Quad1}). In the first part of the paper we'll prove a different, and 
more general (and even shorter) proof of the basic theorem (Theorem 2.5 in \cite{Ca15}) on exact quadratization, showing that, any analytic $\sigma\pi$-ODE is indeed globally exactly quadratizable, and in particular all the solutions $(x,t)\not\in {\rm Dom}\{\Phi\}$ can be calculated as well from the solutions of a quadratic systems, different and smaller dimensional, but still of the same type (\ref{Quad1}).  

Another issue that remained open in \cite{Ca15} is about whether the exact quadratization can be carried over to {\it implicitly defined} ODEs, 
of the type $f(\dot x,x,t)=0$ with both $x$ and $f$ in $\re^n$. For $f$ continuously differentiable with respect to all of its arguments, and having full rank ($n$) with respect to the $\dot x$ indeterminates, the Dini's theorem assures that, at least locally, one can write $\dot x=\phi(x,t)$ for some continuously differentiable $\phi$. We can then talk about a solution $(x,t)$ defined in some domain, in the same sense defined earlier for explicit ODE, carrying it over to implicit ODE as well.   When $f$ is an analytic $\sigma\pi$-function, the local inverse $\phi$ is still not in general 
a $\sigma\pi$-function. Nevertheless, we show in \S \ref{sc:imp} that the local explicit ODE $\dot x=\phi(x,t)$ is $\sigma\pi$-reducible, 
and therefore quadratizable. We perform in this way the extension of the exact quadratization theory, at least up to an extent that includes the basic $\sigma\pi$ case.

Even though it relies on the results of \cite{Ca15}, the present paper is self-contained, and in fact the definition given before for {\it exact quadratization} is new, and, in our view, is more general and suitable for the present purposes\footnote{In \cite{Ca15}, exact quadratization has been defined as a kind of {\it systems immersion}, a concept defined in  \cite{FlKu83} in a framework of control systems}. 
In \S II we review the basic results on exact quadratization, by focusing $\sigma\pi$-ODEs, and introduce some new definitions, among which those of {\it singular, regular} and {\it critical} solution. 
We do not consider more general classes of ODEs, for which we still refer readers to the results of \cite{Ca15}. As a matter of fact, 
the class of ODEs $\sigma\pi$-reducible considered in \cite{Ca15} are all reduced to a $\sigma\pi$-ODE through maps 
$\Phi$ that -- unlike the quadratizing map (\ref{ChV}) , used for $\sigma\pi$-ODEs -- have all ${\rm Dom}\{\Phi\}\subset \cal A$ 
($\cal A$ is the {\it initial} domain, of the ODE to be put, in this case, in $\sigma\pi$-form). In other words, all $\sigma\pi$-reducible ODEs considered in 
\cite{Ca15} are in fact {\it globally} reduced to a $\sigma\pi$-ODE. 

After having shown the globality of the exact quadratization, and its extension to implicit $\sigma\pi$-ODEs, 
which represent the completion of the partial results of \cite{Ca15}, 
in the final part  of this work we come back to the question ({\it q}) we posed at the beginning, and show that the answer is positive. This is the central result of the paper, showing the importance of the exact quadratization theory, 
since a joint application of the latter with it allows 
to write the solution of any analytic $\sigma\pi$-reducible ODE in  its own signature analytic form, a Taylor expansion, 
through an explicit recursive 
calculation of the coefficients from the coefficients of  a quadratization of the original equation. A general formula is provided as well for a lower bound of the convergence ray of the solution Taylor expansion.  

 In the last chapter, we sketch a possible development of the present research as well as show the interplay of the results here presented with 
a differential version of the concept of {\it uniformization} of the 22rd Hilbert's problem.

\section{A review, and Further Results, on the Exact Quadratization for $\sigma\pi$-ODEs}

\subsection{A foreword on terminology: systems, sub-systems, and sub-solutions} 
\label{ssc:ter}

Keeping all definitions and symbols given as yet, hereinafter we talk about a {\it system defined by} an ODE of  the type (\ref{S1}) (hereinafter: 'the ODE') 
as a certain subset (and in particular the whole) of his solutions set. 
If $\cal A$ is the maximal domain of the ODE, a system is given as soon as a subset $\cal S\subset \cal A$ is given, so that we talk about a {\it  system generated by $\cal S$} (associated to the ODE, a specification that often will be omitted) meaning this the set of all solutions $(x,t)\in \cal S$. We name 
the subset $\cal S$ {\it a generating set} of the system. 

We say that $\cal S\subset \cal A$ is an {\it invariance set} 
for the ODE, if  for any 
$(x,t)\in\cal S$, $(x(s),s)\in \cal S$, $\forall s\in {\bf I}_{x,t}$. An invariance set 
$\cal S$ will be also named a {\it domain} (of the ODE).  If $\cal S$ is not an invariance set, 
we can always define the set $\cal S^*\supset \cal S$, said {\it the completion} of $\cal S$, 
collecting all points of solutions graph  passing through a point of $\cal S$. 
Then, by construction, $\cal S^*$ is an invariance set. 
By reason of that, we often identify a system generated by a domain $\cal S$, with the domain 
$\cal S$ itself, and talk about the system as: the {\it system $\cal S$}.  
The largest invariant set  
is $\cal A$, and the ODE itself is a system (the system $\cal A$). Whether 
there isn't any pre-defined $\cal S$, then the system associated to the ODE is understood having domain $\cal A$. 
We also call the system $\cal A$ the {\it total} system. Note that a generating set, or a domain, $\cal S\subset \cal A$ is not necessarily open.
In the following we usually write the ODE component-wise:
\be
\dot x_i=f_i(x,t),\label{S1c}
\ee
understanding that $i\in \cal I=\{1,\ldots, n\}$, where $\cal I$ is named the {\it set of indices} (of the ODE). Let $\cal S\subset \cal A$ be a generating set, and $\cal I^*\subset \cal I$ be a subset of indices. We define the {\it sub-system generated by $(\cal S,\cal I^*)$} 
--- of the system $(\cal A,\cal I)$, which possibly is understood ---
to be the set of all functions $\{x_j(\cdot), j\in \cal I^*\}$ 
 which -- following the definition given in \S I -- are local sub-solutions of  solutions in $\cal S$ with respect to $\cal I^*$  (we'll also refer to 
 them as {\it solutions of} $(\cal S,\cal I^*)$).  If $\cal S$ is a domain, than the sub-system 
generated by $(\cal S,\cal I^*)$ will be simply said: the sub-system $(\cal S,\cal I^*)$. 
Of course, for any given $\cal S_1\subset \cal S_2 \subset \cal A$, system $\cal S_1$ is a subsystem of system $\cal S_2$ and so on, if one 
understands that $\cal S_k$, $k=1,2$, is actually $(\cal S_k, \cal I)$, and $\cal A$ is actually $(\cal A,\cal I)$, in accordance with our definition. 
A sub-system $(\cal S_1,\cal I)$ of $(\cal S_2,\cal I)$ will be also said to be a {\it full} sub-system, whereas if $\cal I^*$ is a proper sub-set of $\cal I$, 
$(\cal S_1,\cal I^*)$ will be said a {\it index-proper} sub-system --- of system $(\cal S_2, \cal I)$, which is also named {\it the hull system} (of system $(\cal S_1,\cal I^*)$). 

Let $\cal S_1,\ldots,\cal S_p\subset \cal S$, and 
$\cal I_1,\ldots,\cal I_q\subset \cal I$ be disjoint partitions of the sets $\cal S$ and $\cal I$ respectively, then 
the collection $\{(\cal S_l,\cal I_m)\}_{l,m}$ will be said: a {\it partition}, and the sub-systems 
generated by $(\cal S_l,\cal I_m)$ will be said {\it parts},   
of the system generated by $(\cal S,\cal I)$.
If $\cal S$ is a domain, we say that the partition 
$\{(\cal S_l,\cal I_m)\}_{l,m}$, is {\it disjoint}, if the sub-systems generated by $({\cal S}_{l_1},{\cal I}_{m_1})$, 
$({\cal S}_{l_2},{\cal I}_{m_2})$, are disjoint for $l_1\ne l_2$,  and $m_1=m_2$, where 'disjoint sub-systems' means that 
they haven't solutions in common. Note that a system partition in {\it full} sub-systems  (namely a {\it full system-partition}) is necessarily disjoint if all the parts are complete (i.e. are defined on sub-domains), but this not necessarily occurs if the parts are index-proper sub-systems. 
Also note that, since every system with generating set $\cal S_l$ has a domain (the completion of $\cal S_l$: $\cal S_l^*$) it follows that 
from any full system-partition one can derive a disjoint partition, namely $\{(\cal S_l^*, \cal I)\}_{l}$, where $l=1,\ldots,p^*\le p$ 
(obviously completing all parts cannot increase the number of parts). 

\subsection{Autonomous sub-systems}
\label{ssc:aut}

A particular kind of  sub-system of a given system $(\cal A,\cal I)$ occurs as the sub-system, say $(\cal S,\cal I^*)$, 
addresses (through $\cal I^*$) equations of the ODE that includes indeterminates $x_j$ with $j\in \cal I^*$ only.  We name these sub-systems {\it autonomous sub-systems}. 
In such a case, the sub-system solutions are actually independent of the components 
$x_j$, for $j\not\in\cal I^*$, of the point $(x,t)\in \cal S$. This means that, denoting:
$$
\cal T=\{t\in \re \quad|\quad \exists x\in \re^n: (x,t)\in \cal A \}
$$
$\cal S$ decomposes as follows
$$
\cal S=\bigcup_{t\in \cal T}\left({\cal P}_t^*\times {\cal Q}_t\times \{t\}\right),
$$
where, if $\iota$ is the cardinality of $\cal I^*$, for any $t\in \cal T$ we have $\cal P_t^*\subset \re^{\iota}$, and ${\cal Q}_t\subset \re^{n-\iota}$.
Thus, provided that the 
canonical injection $\re^{\iota+1}\to \re^{n+1}$ is applied, an autonomous sub-system $(\cal S,\cal I^*)$ (in $\re^{n+1}$) can be 
identified with the system $(\cal P^*, \cal I^*)$ (in $\re^{\iota+1}$), with $\cal P^*=\cup_{t\in \cal T}(\cal P^*_t\times \{t\})$. 

An important case of autonomous subsystem, 
that we encounter in the following, is that given by all sub-solutions zero in ${\cal I}^*$, i.e. by solutions 
of a given system $\cal S$, such that  $x_i\equiv 0$, 
$\forall i\in {\cal I}^*$. We call them 
sub-system {\it projected on ${\cal I}\setminus {\cal I}^*$} which -- provided is non-void --  is described by the scalar ODEs addressed by $\cal I\setminus \cal I^*$, 
modified by setting $x_i=0$ for any $i\in \cal I^*$. Let us denote by $S^0_i$, $i=1,\ldots,n$,  the hyperplane $x_i=0$ in $\re^{n+1}$, and by $S^0$ their union:
\bea
&&S^0_i=\left \{(x,t) \in \re^{n+1}\quad / \quad x_i=0\right\},\label{S0i}\\
&&S^0=\cup_i S^0_i.\label{dfS0}
\eea
The sub-system projected on $\cal I\setminus \cal I^*$ is the index-proper sub-system $(\cal S^*,\cal I\setminus \cal I^*)$, with 
\be
{\cal S}^*=\cal S \bigcap \left( \quad \!\!\bigcup_{j\in \cal I^*}S^0_j\right),\label{barS}
\ee 
of the {\it projected total system} $(\cal A^*,\cal I)$, where $\cal A^*$ is defined by the same identity (\ref{barS}) but the symbol $\cal A$ replaces the symbol 
$\cal S$. Note that $(\cal A^*,\cal I)$ is an index-proper subsystem of $(\cal A,\cal I)$ in $\re^{n+1}$.  
Also, we can consider $(\cal S^*,\cal I\setminus \cal I^*)$ as a {\it full} sub-system of the {\it reduced} system 
$(\cal A^*,\cal I\setminus \cal I^*)$ 
in $\re^{\iota+1}$, as soon as $\cal S^*$ and $\cal A^*$, which are defined as the traces on the hyperplanes $x_j=0$ $j\in 
{\cal I}^*$ of $\cal S$ and $\cal A$ respectively, are identified with subsets of $\re^{n-\iota+1}$ through the canonical injection $\re^{n-\iota+1}\to \re^{n+1}$. 

\subsection{Decomposition of analytic systems in regular and singular parts} 
\label{ssc:sinreg}

\noindent
\begin{definition}\label{df:crs}
{\it For a given system $\cal S$, we name 'critical set'  the set 
\be
{\bf C}={\cal S}\bigcap S^0.\label{Sst}
\ee
The points of ${\bf C}$ are accordingly named: critical points.  
Let us consider the smallest subset of  indices $i_1,\ldots,i_L\in \{1,\ldots,n\}$ such that 
\be
{\bf C}\subset \bigcup_{j=1}^L S^0_{i_j}.\label{senind}
\ee
Then, $i_1,\ldots,i_L$ will be called: criticality indices (of the system). The set of all criticality indices is denoted 
\be
{\cal I}^*=\{i_1,\ldots,i_L\}. \label{ist}
\ee}
\end{definition}

Stated in words, the critical set is the set of point of the ODE domain which are also points of some coordinate hyperplane. 
Along with critical points and criticality indices, we also give the following definition. \\

\noindent
\begin{definition}\label{df:reg} 
{\it A solution, $(x,t)\in \cal S$ 
is said to be 'regular' if the set ${\bf T}_x\subset {\cal T}$ defined as 
\be
{\bf T}_x=\{s\in \cal T\quad /\quad x(s)\in {\bf C}\},\label{CT}
\ee  
has not limit points in $\cal T$. 
}
\end{definition}

We suppose now that the ODE is analytic, i.e. the function $f$ is real-analytic on the domain $\cal A$. 
Thus, the sub-domain $\cal S\subset \cal A$ determines the {\it analytic system} $\cal S$. \\

\noindent
\begin{theorem}\label{th:reg} {\it A solution, $(x,t)$, of an analytic system $\cal S$ 
is regular if and only if his curve passes through a point in 
\be
\cal S_r=\cal S\setminus S^0.\label{Dr}
\ee
}
\end{theorem}

{\it Proof.}  If ${\bf T}_x$ has not limit points in $\cal T$,  then, denoted ${\bf I}_{x,t}$, as usual,  the domain of the solution $(x,t)$,  
we can take an interval $(\alpha,\beta)\subset {\bf I}_{x,t}\subset \cal T$, with 
$(\alpha,\beta)\cap {\bf T}_x=\emptyset$, thus for any $t\in(\alpha,\beta)$, $x(t)\not\in {\bf C}$ and, hence, $x(t)\in \cal S\setminus S^0$.  
Vice-versa, consider a solution passing through $x\in \cal S\setminus S^0$, and let $x(\cdot)$ be his curve . We have 
\be
x_i\ne 0,\quad i=1,\ldots,n.\label{xne0}
\ee
Since the ODE is analytic the curve $x(\cdot)$ is analytic on his domain ${\bf I}_{x,t}$, and 
${\bf T}_x$ is the set of zeroes of the maps $x_i(\cdot)$ (we can consider just the components $i\in {\cal I}^*$).   
All of these maps are not identically zero, by reason of their continuity and by (\ref{xne0}). 
The thesis follows by well known properties of the analytic functions. 
\hfill $\bullet$\\

Theorem \ref{th:reg} implies that any solution of an analytic system either is regular -- hence, it possibly 
crosses the coordinates hyperplanes only at discrete points forming a {\it non dense} sub-set in ${\bf C}$ --  or it lies entirely 
in ${\bf C}$. This motivates the following definition. \\

\noindent
\begin{definition}\label{df:sin} {\it A solution $(x,t)$ of an analytic ODE, is said to be 'singular'  (resp: 'singular in $i$') if 
$x_i\equiv 0$, for some $i\in \cal I^*$ (resp: for $i\in \cal I^*$). The set of all singular (resp: non singular) solutions 
is named 'singular part' (resp: 'regular part') of the $\sigma\pi$-ODE.}
\end{definition}
\vskip0.8cm

Let us denote by $\sigma$ (resp: $\sigma_i$) the set of all singular solutions (resp: in $i$). Then $\sigma$ is the 
singular part of the ODE. Accordingly, we also call $\sigma_i$ the {\it singular part in $i$}, or {\it $i$th singular part}. We have
\be
\sigma=\bigcup_{j\in \cal I^*} \sigma_j.\label{sig}
\ee

Note that, in (\ref{sig}), there well may be $\sigma_j=\emptyset$ (there are not singular solutions in $j$). Now,  
recall that the sub-system projected on $\{i\}$ is the set of all solutions of the ODE obtained by removing the $i$-th equation 
and by setting $x_i=0$ in all the others, thus it is nothing else than the singular part in $i$. That said, 
and accounting of (\ref{barS}), the following proposition ensues.\\

\noindent
\begin{proposition}\label{pr:sin} {\it 
The singular part in $i\in \cal I^*$ of an analytic system $\cal S$ is either empty, or is the sub-system projected on $\{i\}$. The domain of the $i$th 
singular part is
\be
{\bf S}_i=\cal S \bigcap S^0_i.\label{singi}
\ee
}
\end{proposition}
\vskip0.8cm

Let us denote by 
$\cal I_s\subset \cal I$ the set of all indeces $j\in \cal I^*$ such that $\sigma_j\ne \emptyset$. Denote by ${\bf S}$ the domain of $\sigma$, 
by (\ref{sig}) and (\ref{singi}) we have
\be
{\bf S}=\cal S \bigcap \left(\cup_{i\in \cal I_s} S^0_i\right).\label{singdom}
\ee

Since a solution is either regular or singular the following proposition holds.\\

\noindent
\begin{proposition}\label{pr:sreg} {\it  Let be given an analytic system $\cal S$. Define ${\bf R}=\cal S\setminus {\bf S}$. 
Then $\{{\bf S}, {\bf R}\}$ is a disjoint partition of the system $\cal S$, where ${\bf S}$ (resp: ${\bf R}$) is the singular 
(resp: regular) part.}
\end{proposition}
\vskip0.8cm

Thus, by Proposition \ref{pr:sreg} the sub-domain ${\bf R}\subset \cal S$ (the domain of the regular part, or the regular part {\it tout court}, by our usual identification point-solution)  is given, accounting of (\ref{singi}), by
\be
{\bf R}=\cal S \setminus\left(\cup_{i\in  \cal I_s} S^0_i\right)
=\cal S_r \bigcup\left(\cup_{i\in \cal I^*\setminus \cal I_s} S^0_i\right),
\label{regdom}
\ee
where the last identity is obtained by exploiting (\ref{Dr}). Note that the regular part has domain ${\bf R}$, but it can be described 
as well as a system {\it generated by} $\cal S_r$. As a matter of fact, $\cal S_r$ is not (in general) a domain, and ${\bf R}$ is the completion 
of $\cal S_r$. The generating set $\cal S_r$ is a domain if and only if $\cal I^*=\cal I_s$, i.e. all criticality indices are singularity indices, in which 
case we have  $\cal S_r={\bf R}$. The latter case occurs as soon as all solutions $(x,t)\in \cal S_r$  have a curve 
that never crosses a coordinate hyperplane, whereas any solution $(x,t)\in S^0_i$, for some $i$, 
has a curve all lying in $S^0_i$. We denote by $\cal I_{ns}$ the set of all {\it non-singular criticality indices}:
\be
\cal I_{ns}=\cal I^*\setminus \cal I_s.\label{nsci}
\ee
If $\cal I_{ns}\ne \emptyset$, there exist regular solutions $(x,t)\in S^0_i$, for some $i\in \cal I_{ns}$, and we call them {\it critical solutions}.

\subsection{Indices convention}
\label{ssc:ic}

Hereafter we adopt the following convention: {\it a vector} ${\bf v}$ in $\re^\alpha$, $\alpha$ a finite integer, is always a {\it column} vector, and 
${\bf v}'$ his transpose. 
Another convention we make, which will help a lot 
in reducing the length of the most used mathematical expressions, is 
the following {\it indices convention}: if a multi-indexed symbol, say $\xi_{i_1,\ldots,i_n}$, has been defined at some point 
of the paper as a {\it scalar} quantity, then thereinafter the same symbol with the last index removed: $\xi_{i_1,\ldots,i_{n-1}}$ 
shall denote the corresponding {\it column vector} gathering the scalars for all values of $i_n$. Also, if the convention has been applied once, then it can be applied to the aggregated vector again, so that $\xi_{i_1,\ldots,i_{n-2}}$ will be the stack of 
all column vectors for all values of $i_{n-1}$, and so on.

\subsection{Writing $\sigma\pi$-ODEs} 
\label{ssc:wr}
By using this indices convention in (\ref{SP-form1}), we can rewrite a $\sigma\pi$-ODE in the following short form:
\be
\dot x_i=v_i'X_i,\label{SP}
\ee
where $i=1,\ldots, n$. It should be noted that the $i$th equation of a $\sigma\pi$-ODE may be null {\it in many ways}, the extreme possibilities being: if $\nu_i=0$ (which means: there are not monomials), or the vector $v_i=0$, but it cannot be $X_i=0$ though, 
which is never verified. Indeed, $X_{i,l}$ is a {\it function} (of $x$), so with $X_{i,l}=0$ 
we shall mean that $X_{i,l}$ is identically zero (on its maximal domain), which is impossible by definition (there isn't a choice of exponents 
$p^l_{i,j}$ such that $X_{i,l}\equiv 0$ on an open set of $\re^n$).  So, which and how many monomials a $\sigma\pi$-ODE consists of, is a choice 
{\it ad libitum} rather than a material fact. Anyway, in the following, 
as we say that a $\sigma\pi$-ODE is {\it given}, we shall mean that all the monomial numbers $\nu_i$'s, coefficients $v_{i,l}$'s and monomials $X_{i,l}$'s  
have been fixed.

\subsection{Multipartite systems} 
\label{ssc:MPS}

We say that a system in $n$ indeterminates is {\it multipartite}, and, precisely: {\it $m$-partite}, with an integer $m$, if it is double-indexed as follows: 
\be
\dot x_{i,l}=f_{i,l}(x,t),\label{mps}
\ee
where $i=1,\ldots,m$, $l=1,\ldots,\nu_i$ (the index $l$ depends in general of $i$), and $\nu_i$ are integers such that 
$\nu_1+\ldots,\nu_m=n$. Recall that, by the indices convention, $x'=[x_1',\ldots, x_m']$ where $x_i'=[x_{i,1},\ldots,x_{i,\nu_i}]$. 
In accordance with  \S \ref{ssc:ter} we name {\it parts}, the $m$ sub-systems composing (\ref{mps}). In particular, the {\it $i$th part} will denote the sub-system --- 
written in 
{\it vector} form, by applying the indices convention: 
\be
\dot x_i=f_i(x,t).\label{mpsv}
\ee
We say that a system in $n^2$ indeterminates is {\it regularly partitioned}, if it is $n$-partite with all parts  $n$-dimensional.

\subsection{The domain of a $\sigma\pi$-system}

For a  $\sigma\pi$-system the system domain is ${\cal S}={\cal D}\times {\cal T}$, where $\cal T$ is the common 
domain of the coefficients $v_{i,l}$, 
and $\cal D$ is the 
maximal open set on which all the functions $X_{i,l}$, given in (\ref{SP-form1}), 
are well defined and $C^\infty$.  In the following, with some abuse of terminology, we sometimes refer to $\cal D$ (instead of $\cal S$) as the 
system domain. 

Every $\sigma\pi$-system is analytic whereby is $C^\infty$
(and, in particular, locally Lipschitz therein). Consider 
the following (open) sets in $\re^n$, (where $j=1,\ldots,n$): 
\be
{\bf S}^+_j=\{x\in\re^n| x_j>0\};\quad {\bf S}^-_j=\{x\in\re^n| x_j<0\};\label{dfSpm}
\ee
Let us see what is the 
structure of any possible $\cal D$, for any given set of exponents $p_{i,j}^l$.
Let us denote by $ {\cal D}_{i,j}^{l}$ the maximal set where the 
function $\eta_{i,j}^{l}(x)=x_j^{p_{i,j}^{l}}$
is well defined. 
As directly entailed by well known properties of real exponents, there are only four cases:\\

\noindent{\it i)} ${\cal D}_{i,j}^{l}=\re^n$: for any $p_{i,j}^{l}\ge 0$ 
rational, say $p_{i,j}^{l}=n_1/n_2$,  irreducible, and $n_2$ odd;\\
{\it ii)} ${\cal D}_{i,j}^{l}=\overline {\bf S}^+_j$: 
 for any $p_{i,j}^{l}>0$ irrational, or 
rational as in {\it i)} but $n_2$ is even;\\
{\it iii)} ${\cal D}_{i,j}^{l}={\bf S}^+_j$:
as in {\it ii)} but $p_{i,j}^{l}<0$;\\
{\it iv)} ${\cal D}_{i,j}^{l}={\bf S}^+_j\bigcup {\bf S}^-_j$:
 as in {\it i)} but $p_{i,j}^{l}<0$;\\

In \cite{Ca15}, it is shown that, by taking all intersections of the interiors of 
the above sets we obtain $\cal D$, and thus 
we can state the following proposition. \\

\begin{proposition}\label{pr:SPdom} {\it The domain ${\cal D}$ 
of a $\sigma\pi$-system of order $n$, is always an intersection of  a choice of $m\le n^2+1$ sets 
among the following $n^2+1$ sets: $\re^n, {\bf S}^+_j, {\bf S}^+_{j'}\bigcup {\bf S}^-_{j'}$, $1\le j,j'\le n$.}
\end{proposition}
\vskip0.8cm

From Proposition \ref{pr:SPdom} we can see that the domain $\cal D$ of a $\sigma\pi$-system is always an {\it open positive macro-orthant}, 
namely ${\bf O}_{i_1,\ldots,i_L}$, of 
$\re^n$,  which is defined, for a choice $L, i_1,\ldots, i_L\in\{1,\ldots,n\}$, as 
\be
{\bf O}_{i_1,\ldots,i_L}=\{x\in\re^n: x_{i_1} > 0,\ldots , x_{i_L} > 0\},\label{MO}
\ee
to which {\it possibly} the points of some sub-family of the family of the 
{\it coordinate hyperplanes} $x_j=0$ for $j\not\in \{i_1,\ldots,i_L\}$ are drown away. 
 Note that the {\it positive orthant} of $\re^n$ is the macro-orthant ${\bf O}_{1,\ldots,n}$, whereas $\re^n={\bf O}$ (i.e. no choice of  the number $L$). Notice that a $\sigma\pi$-function has always a non-empty domain and all $\sigma\pi$-functions share at least the positive 
orthant in their domains. As a consequence we have the following proposition.\\ 

\begin{proposition}\label{pr:O} {\it Any $\sigma\pi$-ODE has a solution $(x,t)$ for all $x\in {\bf O}_{1,\ldots,n}$ and $t\in \cal T$}
\end{proposition}
\vskip0.8cm

\subsubsection{Example} 
\label{ExdomSP} 

Let us consider the following $\sigma\pi$-ODE with $n=3$, $\nu_1=2$, $\nu_2=1$, $\nu_3=4$ :
\be
\begin{array}{l}
\dot x_1=x_2^{1\over 3}x_3+{2\over \cos(2\pi t)}x_2x_1^{-{1\over 5}}\\
\dot x_2=6x_1x_2^{5}x_3\\
\dot x_3=3x_1^{-8}x_2+4+x_3^{1\over 2}-{3\over 2}x_1x_3
\end{array}\label{eq-exdomSP}
\ee
We have $X_{11}=x_2^{1\over 3}x_3$, $X_{12}=x_2x_1^{-{1\over 5}}$, $X_{21}=x_1x_2^{5}x_3$, 
$X_{31}=x_1^{-8}x_2$, $X_{32}=1$, $X_{33}=x_3^{1\over 2}$, $X_{34}=x_1x_3$, and the coefficients $v_{i,j}$, $i=1,2,3$, $j=1,2,3,4$, are read out at once from (\ref{eq-exdomSP}). The domain of the ODE (\ref{eq-exdomSP}) is $\cal S=\cal D\times \cal T$, with $\cal D$ given, according to the four cases {\it i)-iv)} (cf {\it supra}) :
$$
\cal D=\bigcap_{i=1}^3\bigcap_{l=1}^{\nu_i} {\cal D}^l_{i,j} = \re^3\cap {\cal D}^1_{12} \cap {\cal D}^1_{31} \cap {\cal D}^3_{33}
=\{x\in\re^3\quad | x_1\ne 0, x_3>0\}
$$
(recall that $\cal D^l_{i,j}$ is the domain of $x_j^{p^l_{i,j}}$, i.e. the $j$th term of the $l$th monomial in the $i$th equation, and 
the domains different from $\re^3$ are those of $x_1^{-{1\over 5}}$ in $X_{12}$, $x_1^{-8}$ in  $X_{31}$, and 
$x_3^{1\over 2}$ in $X_{33}$). As for $\cal T$, it is the domain of $v_{12}=2/\cos(2\pi t)$ only: 
$$
\cal T=\bigcup_{k=-\infty}^{+\infty} I_k,\qquad I_k= \{t\in\re \quad | \quad  -\pi+2k\pi < t < \pi + 2k\pi\}
$$
The full system associated to ODE (\ref{eq-exdomSP}), i.e. the set of all solutions, is the union of a countable family of 
sub-systems $\{(\cal D\times I_k, \{1,2,3\})\quad | k:{\rm integer} \}$ associated to the same ODE. In terms of  the set defined in 
(\ref{dfSpm}) and (\ref{S0i}) we have 
$$
\cal D={\bf S}^+_3\setminus S^0_1 \subset \re^3
$$
where the sets (\ref{S0i}) are considered in $\re^n$ (i.e. $\re^3$), neglecting the time coordinate. Thus, $\cal D$ is the macro-orthant of $\re^3$ $x_3>0$, with the points of the hypher-plane $x_1=0$ removed.  
Clearly the sub-systems 
$\cal D\times I_k$ and $\cal D\times I_s$, for $k\ne s$, just differ by translation along the time coordinate. The set ${\bf C}$ defined in (\ref{senind}) (identified with its time-neglecting projection in $\re^3$) is the intersection of $\cal D$ with the coordinate hypherplanes $S^0$ 
thus:
$$
{\bf C}= S^0_2\setminus \{x\in S^0_2\qquad | \quad x_1= 0, \quad \!\!{\rm and}\quad \!\!x_3\le 0\}
$$
and the set of criticality indices is $\cal I^*=\{2\}$. We see that $x_2\equiv 0$ is a sub-solution in the second ODE of (\ref{eq-exdomSP}), thus 
the singular part, ${\bf S}$ is nonempty and equal to the singular part in $\{2\}$, which in turn implies that  ${\bf S}={\bf C}$, 
\footnote{Note that ${\bf S}={\bf C}$ if and only if all criticality indices are singularity indices} and ${\bf R}={\cal D}\setminus {\bf S}$. 

Finally, and coming back to the notation in $\re^3\times I_k$, for any integer $k$ we have that the sub-system 
$(\cal D\times I_k, \{1,2,3\})$ can be decomposed in the two sub-systems $({\bf R}\times I_k, \{1,2,3\})$, and 
$({\bf S}\times I_k, \{1,2,3\})$ (the regular and singular parts). The singular part is the sub-domain given by the half-plane $S^0_2$ for 
$x_3>0$ without the line $x_1=0$, and all solutions of it are given by the projected (on $\{1,3\}$) sub-system: $({\bf S}\times I_k,\{1,3\})$, which is associated 
to the ODE obtained from (\ref{eq-exdomSP}) by setting $x_2=0$ and removing the second equation: 
\be
\begin{array}{l}
\dot x_1=0\\
\dot x_3=4+x_3^{1\over 2}-{3\over 2}x_1x_3
\end{array}\label{ex-singpart}
\ee
The projected system associated to (\ref{ex-singpart}) is regular (it has in fact an empty singular part) \footnote{Notice that the {\it full} system in the indeterminates $x_1,x_3$ associated to 
(\ref{ex-singpart}), has the solution $x_1=0$, and thus a singular part in $x_3$, given by $\dot x_3=4+x_3^{1\over 2}$. However 
$x_1=0$ is not a sub-solution of the sub-system ${\bf S}\times I_k$.} 
and, hence, ${\bf R}\times I_k$ and 
${\bf S}\times I_k$ collect all solutions of the $k$th sub-system $\cal D\times I_k$. 

In conclusion, all solutions $(x,t)$, for $t\in I_k$, and $x_2=0$, $x_1\ne 0$, and $x_3>0$ are singular: they all lie in the half-plane ${\bf S}$. 
All solutions $(x,t)$ with $x_1,x_2\ne 0$ and $x_3>0$ are regular, and lie entirely in ${\bf R}$. Notice that, the regularity domain is composed by 
the four open disjoint ortants of $\re^3$ with $x_3>0$: a solution passing through an $x\in {\bf R}$ remains in the same ortant and never crosses 
the plane $S^0_2$. This is because $S^0_2$ includes the singular part, and all solutions passing through a point of ${\bf S}$ have to lie entirely in ${\bf S}$.\\

\subsubsection{Example} Let us modify ODE (\ref{eq-exdomSP}) by adding the monomial $x_1$ to the second equation. 
Clearly the system domain does not change, and the set of criticality indices is still $\{2\}$. What changes is that the new system has no more a singular part, i.e. it is regular.  In this case solutions $(x,t)$ with $x_2=0$ -- i.e. solutions crossing the hyper-plane  $S^0_2$ -- exist for 
any $x\in \cal D\cap S^0_2$, but none of them lies entirely in $S^0_2$: for an enough small $\epsilon>0$  
the solution $(x,t)\in \cal D\cap S^0_2 \times I_k$ satisfies $x(t\pm \epsilon)\not\in S^0_2$, i.e. {\it escapes} from $S^0_2$, and we can identify 
$(x,t)$ with $(x(t+\epsilon), t+\epsilon)$ (or $(x(t-\epsilon), t-\epsilon)$).\\

\subsection{Exact quadratization of $\sigma\pi$-ODE} 
\label{projODE}

The theorem we present now, is a version of the first statement of Theorem 2.5 in \cite{Ca15}, 
in the formalism of the present paper, and in accordance with the definition of exact quadratization we have given in the Introduction. 
Let be given a $\sigma\pi$-ODE in $n$ indeterminates with $\nu_i$ monomials for $i=1,\ldots,n$, and  
consider the map $\Phi: {\rm Dom}\{\Phi\}\subset \re^{n+1}\to \re^d$ (recall $d=\nu_1+\ldots+\nu_n$):
\be
\Phi_{i,l}(x,t)=x_i^{-1} X_{i,l}.\label{Ph}\
\ee
(actually, constant in the second argument). 
Recall that $l=1,\ldots,\nu_i$, and whether $\nu_i=0$ the corresponding equation is null ($\dot x_i=0$) and there is no $i$-labeled 
entry $\Phi_{i,l}$ in the aggregated vector $\Phi$. Also, recall that the $i$th equation may be zero with $\nu_i\ge 1$, $v_i=0$ 
and $X_{i,l}$ be {\it any} monomial for $l=1,\ldots,\nu_i$  (cf. \S \ref{ssc:wr}). The {\it zero system}, in particular, i.e. the ODE 
$\dot x_i=0$, for $i=1,\ldots,n$, is a $\sigma\pi$-ODE with an arbitrary domain. We conventionally assume $\re^n$ be the domain of the zero system in $n$ indeterminates.\\

\noindent
\begin{theorem} \label{th:EQns}{\it  (Exact Quadratization Theorem -- non singular form) The regular part, 
$\cal S_r=\cal D_r\times \cal T$, of any 
$\sigma\pi$-ODE with domain $\cal S\subset \re^{n+1}$
is exactly quadratizable. 
Moreover, it undergoes the following  quadratization under the map $\Phi^*: (x,t)\mapsto (x,\Phi(x,t))$, with $\Phi$ defined in (\ref{Ph}): 
\bea
&&\dot x_i=(v_i'Z_i) x_i \label{fin}\\
&&\dot Z_{i,l}=\sum_{j=1}^n\pi^l_{i,j} (v'_jZ_j)Z_{i,l} \label{dri}
\eea
where  
\be
\pi^l_{i,j}=p^l_{i,j} -\delta_{i,j},\label{pi}
\ee
$\delta_{i,j}$ being the Kronecker symbol: $=1$  (resp: $=0$) if $i=j$ (resp: $i\ne j$), 
and (\ref{fin}) is (i.e.: remains) the zero equation, if $\nu_i=0$. Finally, the set of all solutions $(\Phi^*(x,t),t)=(x,\Phi(x,t),t)\in \re^{n+d+1}$, with $(x,t)\in {\rm Dom}\{\Phi\}$, of the ODE 
(\ref{fin})--(\ref{dri}) is the system (with domain) $\cal S^*\subset \re^{n+d+1}$:
\be
\begin{array}{rl}
\cal S^*&=\{(x,Z,t)\in \cal D_r\times \re^d\times {\cal T} \quad / \quad Z=\Phi(x,t)\},\\
 &={\cal D}_r\times \Phi({\cal S}_r)\times \cal T=\Phi^*({\cal S}_r)\times \cal T.
 \end{array}
\label{sbdom}
\ee
}
\end{theorem}
\vskip0.8cm

\noindent
\begin{remark}\label{rm:EQ} Of course the ODE (\ref{fin})--(\ref{dri}) is a quadratization according to the definition given in the Introduction, provided 
in the ODE (\ref{Quad1}) we set $m=n+d$, replace in (\ref{fin})--(\ref{dri})  $Z_{i,l}$ with $x_{n+s}$, where 
\be 
s=\alpha_i+l,\quad \alpha_i=\sum_{k=0}^{i-1}\nu_{k} ,\qquad {\rm (}\nu_0=0{\rm )}\label{s}
\ee
and finally replace $v_{i,j}$ of (\ref{Quad1}) with coefficients $v^*_{i,j}$ suitably defined as functions of the parameters $\pi, v$ in 
(\ref{fin})--(\ref{dri}). 
\end{remark}
\vskip0.8cm

{\it Proof of Theorem \ref{th:EQns}.} Let us consider a solution $(x,t)\in {\bf R}$, the domain of the regular part, in accordance with (\ref{regdom}). The analyticity domain 
of the map $\Phi$ in (\ref{Ph}) --- considered as a map $\Phi:\re^{n+1}\to \re^d$ --- is $\cal S_r$ defined in (\ref{Dr}), and by Theorem \ref{th:reg}
we can always choose $(x,t)\in \cal S_r$. Therefore the corresponding curve $s\mapsto x(s)$ is defined and differentiable in an open neighborhood 
of $t$, and the function 
 \be
 Z(t)=\Phi(x(t),t).\label{Zt}
 \ee 
 is defined and differentiable in the same neighborhood. That said, let us calculate the time derivatives of both sides of identity (\ref{Zt}).
 By the chain rule, denoting by $\partial_{z}$ the partial derivative operator with respect to a real indeterminate $z$, and using (\ref{SP}) we have
\be
\dot Z_{i,l}=\sum_{j=1}^n (\partial_{x_j}Z_{i,l})\dot x_j+\partial_t Z_{i,l}=
\sum_{j=1}^n(\partial_{x_j} Z_{i,l})v'_jX_j .\label{Thchr}
\ee
Now, by using (\ref{pi}), and the definition of $\Phi$ given in (\ref{Ph}), we have that $Z$ satisfies, component-wise, the identity:
$
Z_{i,l}=\prod_{j=1}^n x_j^{\pi^l_{i,j}},
$
therefore
$$
\partial_{x_j} Z_{i,l}=\pi_{i,j}^l \left(\prod_{\tau=1}^n x_\tau^{\pi^l_{i,\tau}}\right) x_j^{-1}=
\pi_{i,j}^l Z_{i,l}x_j^{-1},
$$
which yields:
\bea
(\partial_{x_j} Z_{i,l})v'_jX_j &&= \pi_{i,j}^l Z_{i,l}x_j^{-1} (v'_jX_j )=\sum_{s=1}^{\nu_j} \pi_{i,j}^l Z_{i,l}x_j^{-1} v_{j,s}X_{j,s} \nonumber\\
&&=\sum_{s=1}^{\nu_j} \pi_{i,j}^l Z_{i,l} v_{j,s}Z_{j,s} =  \pi_{i,j}^l (v'_jZ_j )Z_{i,l}.\nonumber
\eea
Replacing the expression above in (\ref{Thchr}), the result is (\ref{dri}). The ODE (\ref{fin}) directly comes from the $\sigma\pi$ short form (\ref{SP}) on account of (\ref{Zt}). 
Thus, not only $(x,t)$ is a local sub-solution, in the $x$ components, of the solution $\Phi^*(x,t),t)$ of (\ref{fin})--(\ref{dri}), but 
it holds (\ref{Zt}) as well, at least in a neighborhood of $t$, which proves that $\Phi^*$ is locally compatible.
Therefore all the solutions of (\ref{fin})--(\ref{dri}) 
of the form $(\Phi^*(x,t),t)=(x,\Phi(x,t),t)$ lie in the set $\cal S^*$ defined in (\ref{sbdom}), which is then a sub-domain for the ODE 
(\ref{fin})--(\ref{dri}). This complete the proof of the theorem.
\hfill$\bullet$\\

\noindent
\begin{definition}\label{df:GEQ} {\it We say that an analytic system $\cal S$ 
is 'globally (exactly) quadratizable' if it can be partitioned by a finite collection, $\cal S_1,\ldots,\cal S_L$ of disjoint sub-systems and there exist 
maps $\Phi^1,\ldots,\Phi^L$, with ${\rm Dom}\{\Phi^i\}\cap \cal S_i\ne \emptyset$, such that 
the sub-system $\cal S_i$ is quadratizable under the map 
$\Phi^i$ for $i=1,\ldots,L$. 
}
\end{definition}
\vskip0.8cm

\noindent
\begin{theorem}\label{th:GEQ} {\it  (Exact Quadratization Theorem -- Global form) Any analytic $\sigma\pi$-ODE is globally exactly quadratizable.
}
\end{theorem}
\vskip0.8cm

{\it Proof.} As usual we consider a $\sigma\pi$-ODE in $n$ indeterminates, with a maximal domain $\cal S=\cal D\times \cal T$. 
Moreover we will write ${\bf R}(\cal S), {\bf S}(\cal S)$ for the regular and singular part respectively of a system $\cal S$. 
Let us build up the following sequence of systems: $\cal S_0={\bf S}(\cal S)$, $\cal S_{k+1}={\bf S}(\cal S_k)$ for $k=0,1,\ldots $,  
and let us suppose that $\cal S_k$ is singular (i.e. has a non empty singular part). By 
denoting $\cal I^{\bf s}_{k}$, the singularity indices for the system $\cal S_k$, we have  
$\cal I^{\bf s}_k\ne \emptyset$. 
Now, by definition $\cal S_{k+1}$ is $\cal S_k$ projected on $\cal I\setminus \cal I^{\bf s}_k$, thus  
$\cal S_{k+1}$ is a sub-set of solutions with at least one more zero component than $\cal S_k$, namely $x_i\equiv 0$ with $i\in \cal I^{\bf s}_k$. Therefore there is an $L<+\infty$, such that either $\cal S_L$ is regular (all of its solutions are regular) or it is constituted by the zero 
solution only $x\equiv 0$ (the 'zero system'). By denoting ${\bf R}_k={\bf R}(\cal S_k)$, for $k=0,\ldots, L$,  we have that system $\cal S$ 
 can be partitioned as $({\bf R}_0, \cal S_0)$, where $\cal S_0$ can be partitioned as $({\bf R}_1, \cal S_1)$ and so on, up to 
 $\cal S_{L-1}$ which is partitioned as $({\bf R}_L, \cal S_L)$. Finally,  
 $\cal S$ is partitioned by $({\bf R}_0, {\bf R}_1, \ldots, {\bf R}_{L-1}, \cal S_L)$ where either $\cal S_L=\emptyset$ or it is the zero-system,  which is trivially quadratizable under the identity map\footnote{Recall that conventionally the zero-system has domain $\re^n$.}. The theorem is proven by applying the non-singular form 
of the Exact Quadratization Theorem to the regular systems ${\bf R}_0, {\bf R}_1, \ldots, {\bf R}_{L-1}$. 
\hfill$\bullet$\\

From the proof of Theorem \ref{th:GEQ} we see that any solution non identically zero of an analytic 
$\sigma\pi$-ODE is a regular solution of some 
of its sub-systems, and in particular of  the sub-system projected on the complement of a set of singularity indices. Thus, any non-trivial solution can be 
calculated as analytic prolongation of a sub-solution of a suitable quadratic system. 
For this reason, it will be without loss 
of generality the assumption we make hereinafter: 
as we talk about 'a solution' of an analytic $\sigma\pi$-ODE, we always understand, unless explicitly specified, that it is a 
{\it regular} solution.\\

\subsection{Exact quadratization for implicitly defined $\sigma\pi$-systems} 
\label{sc:imp}
Let us consider the {\it implicitly defined} ODE:
\be
f(x,\dot x,t)=0,\label{IDO}
\ee
where $x\in \re^n$ and $f: U {\rm (open)}\subset \re^{2n+1}\to \re^n$ is a continuously differentiable function having full rank ($n$)
in the $\dot x$ indeterminates, at a point $z=(x,\dot x,t)\in U$, i.e. $(\partial_{\dot x} f)(x,\dot x,t)$ has rank $n$ (with 
$\partial_{\dot x}$ denoting the Jacobian with respect to $\dot x$). Let $\cal U$ be the set in $U$ where the algebraic equation 
(\ref{IDO}) is satisfied:
\be
\cal U=\{(x,\dot x, t)\in U\quad | \quad f(x,\dot x,t)=0\}.\label{calU}
\ee
Then, 
by the Dini's theorem (in many indeterminates), for any $(x,\dot x, t)\in \cal U$
there exist two  
open sets $U_1,U_2$,   $(x,\dot x,t)\ni U_1\subset U$, $(x,t)\ni U_2\in \re^{n+1}$, 
and a continuously differentiable 
function $\phi:U_2 \to \re^n$, such that, denoting $y=(x,t)$: 
\bea
&&\cal U\cap U_1=\{(x,\phi(x,t),t) \quad | \quad (x,t)\in U_2 \},\label{If1}\\
&&(\partial_y\phi)(x,t) =-\Big((\partial_{\dot x} f)(x,\phi(x,t) ,t)\Big)^{-1}\left((\partial_{y} f)(x,\phi(x,t) ,t)\right), \qquad \forall y\in U_2.\label{If2}
\eea
Therefore, the {\it explicit} ODE 
\be
\dot x=\phi(x,t),\label{expO}
\ee
on the domain $U_2\subset \re^{n+1}$ has solutions that are local solutions of the implicit ODE (\ref{IDO}). 
If $(\partial_{\dot x} f)(z)$ has full rank $\forall z\in U$, then all solutions of (\ref{IDO}) are local solutions of some explicit ODE of the type 
(\ref{expO}).  Then we can naturally extend 
the definition given earlier of {\it system} (associate to an explicit ODE) to {\it implicit systems} (associated to an implicit ODE),  
by defining an implicit system as a {\it domain} 
$\cal S\subset \re^{n+1}$, and a solution $(x,t)\in \cal S$, as the maximal map $s\mapsto x(s)$, with $x(t)=x$, 
which solves the implicit ODE (\ref{IDO}). From what said earlier it is $\cal S\ne \emptyset$, since for any $(x,\dot x,t)\in \cal U$ there is 
always an explicit ODE  locally solved by $(x,t)$. 
 
These facts carry over to the analytic case, and we have an analytic local inverse $\phi$ associated an analytic implicit system function $f$ at each point of its own domain. Now, the question we pose is the following: suppose $f$ is an analytic $\sigma\pi$-function 
in $2n$ indeterminates, 
which means it has the well known $\sigma\pi$ expression with respect to the indeterminates $x,\dot x$, for some set of analytic coefficients 
$v_{i,l}$, then is any associated local inverse $\phi$ an analytic $\sigma\pi$-function as well? The answer is in general negative, but the following Theorem shows that the explicit local system given by (\ref{expO}) 
is always $\sigma\pi$-{\it reducible} (cf. \S I) which is sufficient for our present purposes.\\ 

\noindent
\begin{theorem} \label{th:imp}
{\it  Let be given an analytic $\sigma\pi$-ODE in the implicit form (\ref{IDO}), with $f$ of full rank ($n$) with respect to the 
$\dot x$ indeterminates, at a point $z=(x,\dot x,t)\in \cal U$.  Then, the associated explicit local representation (\ref{expO}) is $\sigma\pi$-reducible. 
}
\end{theorem}
\vskip0.8cm

Before proving the Theorem, some comments are in order. Recall that with '$\sigma\pi$-reducible' we mean that 
the ODE (\ref{expO}) has solutions that are all sub-solutions of a (larger-dimensional) $\sigma\pi$-ODE, i.e. there are additional 
indeterminates  $x_{n+1},\ldots,x_{n+d}$, for some $d\in\na$, and analytic time-functions $v_{i,l}$ such that: 
\be
\dot x_i= \sum_{l=1}^{\nu_i}v_{i,l}\prod_{j=1}^{n+d} x_j^{p^l_{i,j}},\qquad i=1,\ldots,n+d,\label{SPes}
\ee
such that any solution $(x,t)\in \cal D\subset \re^{n+1}$ where $\cal D$ is the domain of (\ref{expO}) is a local sub-solution 
of (\ref{SPes}) in the first $n$ indeterminates. We can also say that the system $\cal D$ associated to ODE (\ref{expO}) is 
$\sigma\pi$-reducible if it is the sub-system $(\cal D,\{1,\ldots,n\})$ associated to the ODE (\ref{SPes}). Clearly, a 
$\sigma\pi$-reducible system undergoes a quadratization, and, hence, by Theorem 2 we can naturally 
build up the extension of the exact quadratization theory to implicit analytic full-rank $\sigma\pi$-systems of the type (\ref{IDO}). \\

Some classes of $\sigma\pi$-reducible (explicit) systems have been studied in \cite{Ca15}, these classes of systems 
are defined by the class of functions the system-function belongs to.  We here recall one of them, which is very basic, 
referring reader to \cite{Ca15} for more detail: (system function given by) a  
{\it finite composition} of $\sigma\pi$-functions, which in \cite{Ca15} is named, with some abuse of terminology, 
an {\it algebraic} function. For instance, consider the following ODE: 
\be
\begin{array}{l}
\dot x_1= {x_1^2x_2^{-{1\over 4}} \over x_1+x_2^2}\\
\dot x_2=x_1^2x_2
\end{array}\label{exNSP}
\ee
the first equation is not $\sigma\pi$, but the system function can be written as $\phi_1\circ \phi_2$, 
where  $\phi_1(x_1,x_2,x_3)=x_1^2x_2^{-{1\over 4}}x_3^{-1}$ 
and $\phi_2(x_1,x_2)=(x_1,x_2, x_1+x_2)$, and thus, it is the composition of two  $\sigma\pi$-functions.  By taking the time derivative of $x_3=x_1+x_2^2$ we have $\dot x_3=\dot x_1+2x_2\dot x_2$, 
which is $\sigma\pi$, since $\dot x_1=\phi_1(x_1,x_2,x_3)$ is $\sigma\pi$, $\dot x_2$ is $\sigma\pi$ by the second of (\ref{exNSP}), 
and sums and products of $\sigma\pi$-functions
results in a $\sigma\pi$-function. By adjoining the equation for $x_3$ to (\ref{exNSP}), and replacing $x_3$ for $x_1+x_2^2$, we get 
a $\sigma\pi$-ODE in the three indeterminates $x_1,x_2,x_3$. We recall also this (quite straightforward) result: let $\phi$ a scalar $\sigma\pi$-function in $x\in\re^n$, then $\partial _{x_1}f$, the partial derivative of $f$ in $x_i$, is a $\sigma\pi$-function (in $x$). \\ 

{\it Proof of Theorem \ref{th:imp}.} By the Dini's theorem we get the explicit ODE (\ref{expO}) whose domain $U_2$ is such that every 
$(x,t)\in U_2$ is a local solution of (\ref{IDO}) as well. Moreover $\phi$ is continuously differentiable, thus we can set $z=\phi(x,t)$ and 
take the time-derivative, the result is
\be
\begin{array}{l}
\dot x=z,\\
\dot z=(\partial_x\phi) \cdot \dot x+\partial_t \phi =(\partial_x\phi) \cdot z+\partial_t\phi
\end{array}\label{tran}
\ee
In (\ref{If2}), note that, for $y=(x,t)$, $\partial_y\gamma = (\partial_x \gamma,  \partial_t\gamma)$ (where the symbol $\gamma$ stands for 
$\phi$ or $f$) thus (\ref{If2}) can be rewritten as
\be
\begin{array}{l}
(\partial_x\phi)(x,t) =-\Big((\partial_{\dot x} f)(x,\phi(x,t) ,t)\Big)^{-1}\left((\partial_{x} f)(x,\phi(x,t) ,t)\right), \\
(\partial_t\phi)(x,t) =-\Big((\partial_{\dot x} f)(x,\phi(x,t) ,t)\Big)^{-1}\left((\partial_{t} f)(x,\phi(x,t) ,t)\right).
\end{array}
\label{If2}
\ee
Now, since $f(x,\dot x,t)$ is $\sigma\pi$ with respect to $x$ and $\dot x$, 
any component of the matrix $\partial_{\dot x}f$ or $\partial_{x}f$ is a $\sigma\pi$-function in $x,\dot x$ as well. 
Moreover, any component of $(\partial_xf)^{-1}$ is the ratio of the determinant of some minor and the determinant 
of the whole matrix, thus a ratio of polynomials -- hence a $\sigma\pi$-function -- 
of the components of $\partial_xf$. Therefore, we can set 
\be
-\Big((\partial_{\dot x} f)(x,\phi(x,t) ,t)\Big)^{-1}\left((\partial_{x} f)(x,\phi(x,t) ,t)\right)=h(x,\phi(x,t),t), \label{conc1}
\ee 
where the matrix function $h$ has  finite compositions of $\sigma\pi$-functions in all entries. The function $\partial_tf$ has the usual $\sigma\pi$ expression 
(\ref{SP-form1})  
but the coefficients $\dot v_{i,j}$ replaces $v_{i,j}$, therefore is an analytic $\sigma\pi$-function, and we can set 
\be
-\Big((\partial_{\dot x} f)(x,\phi(x,t) ,t)\Big)^{-1}\left((\partial_{t} f)(x,\phi(x,t) ,t)\right)=g(x,\phi(x,t),t)., \label{conc2}
\ee 
and $g$ is still a finite composition of $\sigma\pi$-functions. By (\ref{If2})--(\ref{conc2}) eq.s (\ref{tran}) rewrite 
\be
\begin{array}{l}
\dot x=z,\\
\dot z=h(x,z,t)  z+g(x,z,t)=\psi(x,z,t)
\end{array}\label{tran2}
\ee 
 where $\psi$ is clearly a finite composition of $\sigma\pi$-functions. In conlusion, the ODE (\ref{tran2}) is $\sigma\pi$-reducible, and since 
 the solution $(x,t)$ of the explicit ODE (\ref{expO}) is a sub-solution of (\ref{tran2}), we conclude that (\ref{expO}) is $\sigma\pi$-reducible 
 as well. \hfill $\bullet$\\
 
 Theorem \ref{th:imp} implies in turn the following corollary. \\
 
 \begin{corollary}\label{co:imp}
 {\it Any implicit analytic $\sigma\pi$-ODE of the type (\ref{IDO}), $f$ having rank $n$ in the $\dot x$ indeterminates, 
 for any $z=(x,\dot x,t)\in \cal U$, is globally quadratizable.
 }
\end{corollary}
\vskip0.8cm
 
 It should be stressed that,  unlike the case of an explicit $\sigma\pi$-ODE,  
 any solution $(x,t)$ in the domain of an implicit analytic $\sigma\pi$-ODE is, by Corollary \ref{co:imp},  
 a local sub-solution of a quadratic system of the type (\ref{Quad1}), but the latter is not in general unique,  and in fact depends of the solution 
 $(x,t)$ itself. \\
  
 \begin{remark}
 \label{rm:impl}
 A more important remark is that, even if the function $\phi$ of the locally defined explicit ODE (\ref{expO}) could be difficult to 
 calculate, we have not to perform this calculation though in order to calculate a solution of the original implicit ODE. 
For,  by the proof of Theorem \ref{th:imp} we realize that any ODE solution is a sub-solution (in the first $n$ components) 
 of the explicit ODE (\ref{tran2}), which is defined by a vector field that can be calculated through (\ref{conc1}), and (\ref{conc2}), from the initial 
 implicit vector field $f$: we can then quadratize system 
 (\ref{conc1}), (\ref{conc2}) and the solution of the original implicit system is {\it a fortiori} 
 a local sub-solution of the quadratized ODE, thus still recoverable by analytic extension.
 \end{remark}
 \vskip0.8cm

\subsection{The structure of the canonical quadratization}
\label{sc:SCO}
As for the ODEs (\ref{fin}) and (\ref{dri}), we here follow the terminology introduced in \cite{Ca15},  but 
accounting as well of the new definitions given in the previous sections. 
Hereinafter we name the ODE (\ref{SP}), with domain $\cal S$, 
the {\it original $\sigma\pi$-system}.\\ 

\begin{definition}\label{df:quahu}{\it We name the system (on $\re^{n+d}\times {\cal T}$) defined by the ODEs (\ref{fin}) and (\ref{dri}), 
the quadratic hull (of the original $\sigma\pi$-system)}
\end{definition}
 \vskip0.8cm

\begin{definition}\label{df:quaim} {\it We name 'quadratic image' or 'image system' or 'canonical quadratization' 
(of a given $\sigma\pi$-system, that is referred to as 'original system') -- the system defined by the ODEs 
(\ref{fin}) and (\ref{dri}), on the sub-domain ${\cal S}^*$ defined by (\ref{sbdom}). 
Moreover, the two sub-systems of the quadratic image, given by the ODEs (\ref{fin}) and (\ref{dri}), are  
said the 'final stage' (of the quadratization) and the 'Driver' (of the original system) respectively.
}\footnote{We chose, in the paper \cite{Ca15}, the names {\it Driver} and {\it final stage}, since a canonical quadratization can be 
depicted, in a similarity with electrical devices, as a {\it cascade} of systems, where the Driver actually 'drives' the final stage yielding the signal representing the original ODE solution.}
\end{definition}
 \vskip0.8cm

Note that, 
in accordance with \S \ref{ssc:MPS} the Driver is a {\it $n$-partite} system, whose $i$th part consists of $\nu_i$ (scalar) equations. 
Moreover is an {\it autonomous} system, i.e. the sub-system $({\cal S}^*, {\cal I}^{\bf d})$ of the image system, where $\cal I^{\rm d}\subset \cal I$ the sub-set of 
indices addressing the last $d$ equations of the image system, thus following \S \ref{ssc:aut}, and accounting as well of the structure of 
$\cal S^*$ given in (\ref{sbdom}), it can be identified 
with the system (in $\re^d\times \cal T$): $(\Phi(\cal D_r\times \cal T)\times \cal T, \cal I^{\rm d})$. We can omit the (obvious) symbol $\cal I^{\rm d}$, and consider the Driver with the full system: $\Phi(\cal D_r\times \cal T)\times \cal T$. 

In the following, as we talk about a Driver, it shall be always entailed {\it of which} original ($\sigma\pi$) system, that Driver is. 
By reason of that, we also name {\it pre-Driver} the original $\sigma\pi$-system, insofar we move the focus from the original system to 
the associated Driver. In order to refer ODE (\ref{dri}) as a sub-system of the quadratic hull we use the different term of {\it Driver ODE}, instead of the single word 'Driver'. 

The canonic quadratization holds under the map $\Phi^*$, which we call the {\it canonic quadratizing map} or 
the {\it quadratic immersion}.   

\subsection{Non-canonical quadratizations, and some examples.}
\label{sc:SCO}

According to the definition given in the previous paragraph, the Driver is a part of the canonical quadratization, as well as the final stage.  
Nevertheless,  there is a kind of quadratization, that can be always applied to the same original system, in which the final system is redundant, and can be considered consisting of the sole Driver.  \\

\begin{definition}\label{df:inc-qua} {\it We say that the canonical quadratization (of a given $\sigma\pi$-ODE) is 'inclusive', if all solutions of the original ODE are local sub-solutions of the Driver.}\end{definition}
 \vskip0.8cm

By the definition given above, a canonical quadratization is inclusive if and only if for $i=1,\ldots,n$,  $\Phi_{i,l}=x_i$ for some $l=1,\ldots, \nu_i$, which is verified if for all $i$ there is an $l$ such that $X_{i,l}=x_i^2$. If for some $i$ the original equation has not the monomial $x_i^2$ 
one can always add the fictitious term $0\cdot x_i^2$ in the $i$th equation. 
The addition of a fictitious term obviously doesn't change the solution of the ODE (\ref{SP-form1}), but in the canonical quadratization it does change the  dimension ($d$) of the target space for the aggregate map $\Phi$ given by (\ref{ChV}), 
and in fact changes the map $\Phi$ itself, making redundant the quadratizing map $\Phi^*$. Indeed, we can easily verify, by applying formula (\ref{dri}) that the driver equation for 
$Z_{i,l}=\Phi_{i,l}=x_i$ (thus, associated to the monomial $X_{i,l}=x_i^2$) is equal to the $i$th equation of the final sub-system (\ref{fin}). 
For, let us add $0\cdot x_i^2$ at the $i$th equation as the last term, then the new $\sigma\pi$-ODE has $\nu_i^*=\nu_i+1$ monomials, and 
the Driver has the further variable $Z_{i,\nu_i^*}=x_i^{-1}x_i^2=x_i$. By (\ref{dri}) we have
$$
\dot Z_{i,\nu_i^*}=\sum_{j=1}^n \pi^{\nu^*_i}_{i,j} (v_j'Z_j) Z_{i,\nu^*_i},
$$
but  from (\ref{pi}) $\pi^{\nu^*_i}_{i,j}=\delta_{i,j}$, ($X_{i,\nu^*_i}=x_i^2=x_1^0\cdots x_{i-1}^0x_i^2x_{i+1}^0\cdots x_n^0$ from which $p_{i,j}=0$ for $i\ne j$ and $p_{ii}=2$) and thus $\dot Z_{i,\nu_i^*}=(v_i'Z_i)Z_{i,\nu_i^*}$ which is the $i$th final equation (\ref{fin}), since 
$Z_{i,\nu_i^*}=x_i$ and accounting as well that, with respect to the canonical quadratization, 
$v_i$ has the further last component $0$ and $Z_i$ has the further component $Z_{i,\nu_i^*}$, which leaves the product 
$v_i'Z_i$ unchanged. 

In conclusion, given any $\sigma\pi$-ODE, we naturally obtain an inclusive quadratization by applying  the canonical quadratization to the ODE modified by adding, for $i=1,\ldots,n$, the fictitious term $0\cdot x_i^2$ in the $i$th equation, provided the monomial $x_i^2$ is not therein yet present.   
Clearly, the new quadratization has as many new Driver variables, say $Z_{i,\nu_i+1}=x_i$, as the original ODE has equations $i$-labeled lacking the monomial $x_i^2$. Correspondingly, a subset of the final equations in the canonical quadratization turns out to be 
replied among the new Driver equations.  Since the natural inclusive quadratization is {\it redundant}, we usually omit the final stage, thus for an inclusive quadratization we assume that the quadratic hull 
of the original $\sigma\pi$-ODE is the Driver ODE, and identify the quadratization with the sole Driver.\\

\subsubsection{Example}
\label{Exnoncan}
Let us consider the following $\sigma\pi$-ODE with $n=3$, $\nu_1=2$, $\nu_2=1$, $\nu_3=2$: 
\bea
&&\dot x_1=x_2x_3+2x_1^{-1/3}\nonumber\\
&&\dot x_2=4x_1x_2^4x_3,\nonumber\\
&&\dot x_3=5x_1^{-3} x_2+3.\nonumber
\eea
The canonical quadratization is composed by a Driver in 5 indeterminates (as many as is the number of  monomials) defined as 
\bea
&&Z_{11}=x_2x_3x^{-1}_1; \quad Z_{12}=x_1^{-4/3},\nonumber\\
&&Z_{21}=x_2^5x_3,\nonumber\\
&&Z_{31}=x_1^{-3}x_2x_3^{-1}; \quad Z_{32}=x_3^{-1}.\nonumber
\eea
By (\ref{fin} ), the final ODE is: 
\be
\begin{array}{ll}
\dot x_1=(Z_{11}+2Z_{12})x_1 &{\rm (} =(v_1'Z_1)x_1 {\rm )}\\
\dot x_2=4Z_{21}x_2, &{\rm (} =(v_2'Z_2)x_2 {\rm )}\\
\dot x_3=(5Z_{31}+3Z_{32})x_3, &{\rm (} =(v_3'Z_3)x_3 {\rm )}.\\
\end{array}\label{fex}
\ee
By (\ref{dri}) we calculate the Driver ODE. A trick for calculating the equation of $Z_{i,l}$ is looking at the monomial which defines it: 
just write the linear composition of $v_j'Z_j$ using the corresponding powers of $x_j$ in the monomial as coefficients, and finally multiply by $Z_{i,l}$: 
\bea
&&\dot Z_{11}=(v'_2Z_2+v'_3Z_3-v_1'Z_1)Z_{11}; \quad \dot Z_{12}-(4/3)(v_1'Z_1)Z_{12},\nonumber\\
&&\dot Z_{21}=(5v_2'Z_2+v_3'Z_3) Z_{21},\nonumber\\
&&\dot Z_{31}=(-3v_1'Z_1+v'_2Z_2-v_3Z_3)Z_{31}; \quad \dot Z_{32}=-(v_3'Z_3)Z_{32},\nonumber
\eea
and the $v_j'Z_j$, $j=1,2,3$, are given from (\ref{fex}). The natural inclusive quadratization is obtained by simply introducing the 
new Driver indeterminates: $Z_{13}=x_1$, $Z_{22}=x_2$, $Z_{33}=x_3$, whose equations are directly given by (\ref{fex}) 
upon substituting $x_1,x_2,x_3$. \\

Adding a fictitious monomial is a method that was not investigated in \cite{Ca15}: we now show that certain side results of \cite{Ca15} can be 
derived in a straightforward way, and even generalized, by using this method. Let us consider a general monomial function of the $n$ indeterminates of a given 
$\sigma\pi$-ODE, namely $g(x)\in \re$: 
\be
g(x)=x^{q_1}_1x^{q_2}_2\cdots x^{q_n}_n;\qquad q_1,\ldots,q_n\in \re,\label{funcg}
\ee
By using the exact quadratization we can prove that $g$ satisfies a scalar quadratic ODE: let us add the monomial $0\cdot g$ in the $i$th equation, and apply the canonical quadratization. With the same symbols as above, we have $Z_{i,\nu_i^*}=x_i^{-1}g$, that is 
$g=Z_{i,\nu_i^*}x_i$, and by taking the time derivative:
$$ 
\dot g=\dot Z_{i,\nu_i^*}x_i+Z_{i,\nu_i^*}\dot x_i
$$
and by the Driver and final equations we have 
\be
\begin{array}{rl}
\dot g&=\left(\sum_{i\ne j=1}^n q_j (v_j'Z_j) +(q_i-1) v_i'Z_i\right) Z_{i,\nu^*_i}x_i+ (v_i'Z_i)x_iZ_{i,\nu^*_i}\\
 &=\left(\sum_{j=1}^n q_j (v_j'Z_j)\right) g=h\cdot g
 \end{array}\label{eqg}
\ee
Let $(x,t_0)\in \cal D_r\times \cal T$ be a solution of the $\sigma\pi$-ODE, and $I_{x,t_0}\subset \cal T$ be its maximal open interval of definition. Then, there is an interval $t_o\ni I$(open) $\subset I_{x,t_0}$ where the time-functions 
$$
h_t= \sum_{j=1}^n q_j (v_j'(t)Z_j(t)); \qquad g_t=g(x(t))
$$  
are well defined (and, hence, analytic). Using (\ref{eqg}), $h_t$ in $I$ is given by:
\be
g_t=e^{\int_{t_0}^t h_\tau d\tau} g_{t_0},\label{gsolexp}
\ee
which allows the calculation of $g_t$ from a Driver solution $Z(t)$ (through $h_t$). If $g(x)=x_i$, 
and thus $q_j=0$ for $j\ne i$ and $q_i=1$, (\ref{eqg}) turns out into the well known $i$th final equation (we come back in this way to the canonical quadratization). In the latter case, $g_t=x_i(t)$, and it is defined by definition on the whole of $I_{x,t_0}$, which is in general a larger interval than $I$.  Thus, the exponential solution formula (\ref{gsolexp}) gives $x_i$ only on the interval $I$ (where $Z(t)$ is defined). Outside of $I$, $x(t)$ is determined though, by analytic prolongation. 

We have found out in a very simpler way the same quadratization shown in \cite{Ca15}, where a further system: the {\it medial system} 
given by the ODE (\ref{eqg}) for $g=X_{i,l}$, is  'driven' by $Z(t)$ (by the Driver) and gives all the monomial functions $X_{i,l}(t)$, which in turn give the original solution $x(t)$ through integral action:
$$
x_i(t)=\int_{t_0}^t v_i'(\tau)X_i(\tau) d\tau
$$
We name this kind of quadratization: a {\it medial quadratization}. \\

\subsubsection{Example}
\label{Ex2}
Let $a,b$ be two fixed positive real numbers, $t_0\in\re$, and let us consider the affine ODE:
\be
\dot x=-a x-b.\label{ex-lsys}
\ee
In our terminology, we are considering the system $(\re^2,\{1\})$ associated to ODE (\ref{ex-lsys}). Since $b\ne 0$ there aren't singular solutions, and the regularity domain is $\cal S_r=\cal D_r\times \re$ with 
$\cal D_r=\re^+\cup \re^-$.  
As a $\sigma\pi$-system, (\ref{ex-lsys}) is characterized by $n=1$, $\nu_i\equiv \nu_1=2$, $v_{1,1}=-a$, $v_{1,2}=-b$, 
$X_{1,1}=x$, $X_{1,2}=1$, $p^1_{1,1}=1$, $p^2_{1,1}=0$. The components of the canonical Driver are 
$Z_{1,1}=1$, $Z_{1,2}=x^{-1}$, and by (\ref{dri}) we obtain the canonical Driver equations:
\bea
&&\dot Z_{1,1}=0,\label{ex-dri1}\\
&&\dot Z_{1,2}=aZ_{1,1}Z_{1,2}+bZ^2_{1,2}.\label{ex-dri2}
\eea
By (\ref{fin}) we get the final equation:
\be
\dot x=(-aZ_{1,1}-bZ_{1,2})x.\label{ex-fin}
\ee
The Driver map is $\Phi(x,t)=(1,x^{-1})$, which has domain $\cal S_r$.  The quadratic hull is the system, with domain $\re^4$, associated 
to the three equations (\ref{ex-dri1})-(\ref{ex-fin}), whereas the image system is given by the same ODEs  on the 
domain $\re\times \Phi(\cal S_r) \times \re$.

Let us consider the solution $(x,t_0)$ of (\ref{ex-lsys}), and denote $t\mapsto x(t)$ the corresponding 
solution curve. We can easily calculate $x(t)$: 
\be
x(t)=x e^{-a(t-t_0)}-{b\over a} (1-e^{-a(t-t_0)}),\label{solex}
\ee
and we have ${\bf I}_{x,t_0}=\re$. Let us define 
${\bf T}_{x,t_0}=\{t\in\re: x(t)=0\} $ 
which is the set (\ref{CT}) for the case at issue. ${\bf T}_{x,t_0}$ is a countable set, by Theorem \ref{th:reg}, and we can calculate it 
by distinguishing three cases:\\

{\it (i)}  $x\ge 0$. In this case ${\bf T}_{t_0,\bar x}=\{\bar t\}$, with 
\be
\bar t=t_0+{c\over a};\qquad c=-\ln\left({b\over a\bar x+b}\right) >0.\label{ex-tmp}
\ee
and $t_0\le \bar t$.\\

{\it (ii)}  $-(b/a)<x<0$. In this case ${\bf T}_{t_0,\bar x}=\{\bar t\}$, with $\bar t$ as 
in (\ref{ex-tmp}), but $t_0>\bar t$.\\

{\it (iii)} $x\le-(b/a)$. In this case 
${\bf T}_{0,\bar x}=\emptyset$.  \\  

Now, consider the solution $(\Phi(x,t),t)=(1,x^{-1},t)$ of the Driver. 
By (\ref{ex-dri1}) we have $Z_{1,1}\equiv 1$ for any $t$. The component $Z_{2,2}$ is given by  
the ODE (\ref{ex-dri2}), which  is a {\it Bernoulli} quadratic differential equation, and thus the general solution can be calculated by a well known method: define $W_{1,2}=Z^{-1}_{1,2}$, then (as it is easy to verify) $W_{1,2}$ satisfies the {\it affine} equation:
\be
\dot W_{1,2}=-aZ_{1,1} W_{1,2}-b.\label{ex-W12}
\ee
Notice that, since $Z_{1,1}=1$, (\ref{ex-W12}) is the same equation as the original system (\ref{ex-lsys}). 
Moreover $W_{1,2}(t_0)=Z^{-1}_{1,2}(t_0)=x$.  
The solution is then equal to (\ref{solex}), and by $Z_{1,2}=W^{-1}_{1,2}$, we 
get the general expression for $Z_{1,2}$:
\be
Z_{1,2}(t)={a \over  (ax +b)e^{-a(t-t_0)}-b}.\label{ex-SG}
\ee
The above function is (as expected) not defined for $t\in {\bf T}_{x,t_0}$. 
We can use formula (\ref{gsolexp})  with $g(x)=x$ in order to express the solution given by the final system, i.e. ODE (\ref{ex-fin}) 
{\it driven} by the Driver: 
\be
x(t)=e^{\int_{t_0}^t Z_{1}^Tv_{1}} x=e^{-\int_{t_0}^t (aZ_{1,1}+bZ_{1,2})} x,
\label{ex-solfinal1}
\ee
If $x\le -(b/a)$ -- case {\it (iii)} above --  (\ref{ex-SG}) is always defined and non-zero 
on $\re$  and (\ref{ex-solfinal1}) holds $\forall t\in \re$, and gives back the time-function (\ref{solex}) $\forall t$. Thus. according to the Exact Quadratization Theorem, $x(t)$ is a sub-solution of the quadratic image, and for $x\le -(b/a)$ in particular is a {\it global} sub-solution (cf \S I) 
defined on ${\bf I}_{x,t_0}=\re$.  If $x\ge 0$ (resp: $-(b/a)<x<0$) --  case {\it (i)} (resp: {\it (ii)})  -- the solution $(x,t_0)$ of the original ODE is a  local sub-solution defined on the 
sub-interval ${\bf I}=(-\infty, \bar t )$,  (resp: $(\bar t, +\infty)$) with $\bar t$ given by (\ref{ex-tmp}), and 
(\ref{ex-solfinal1})  -- calculated for $t\in {\bf I}$ -- gives the restriction on ${\bf I}$ of the function (\ref{solex}). 
It should be stressed that, a strict-sense local sub-solution of the original solution $(x,t_0)$ -- which occurs as $x>-(b/a)$ -- is given by 
(\ref{ex-solfinal1})  only in ${\bf I}$, and it couldn't be otherwise since the time function defined by (\ref{ex-solfinal1}) is either positive or negative (which depends of the sign of $x$) everywhere in ${\bf I}$. Outside ${\bf I}$ though, it shall be still equal to $(x,t_0)$ by analytic prolongation.\\

\subsection{Inverse Driver}
\label{Inv-driv}

By the Exact Quadratization Theorem (any of the two forms) a Driver solution is a {\it $Z$-sub-solution} of the quadratization 
i.e. the components from 
the $(n+1)$th to the $(n+d)$th of a solution lying in the domain $\cal S^*$ defined in (\ref{sbdom}). Thus, we have pointwise 
$Z(t)=\Phi(x(t))$ for any $t\in I_{x,t}$ such that $x(t)\in \cal D_r$, and hence --- by (\ref{Ph}) ---  $Z^{-1}(t)$ is defined for the same 
set of $t$ as $Z(t)$ i.e. $\forall t\in I_{x,t}$ with the possible exception of a subset that hasn't limit points.  
in $I_{x,t}$. Therefore, we can define 
\be
W_{i,l}=Z_{i,l}^{-1}\quad,\label{dfW}
\ee
so everywhere $Z$ is defined and differentiable so does $W$ as well. By taking the time derivative, 
we have that $W_{i,l}$ satisfies the equation:
\be
\dot W_{i,l}=-Z^{-2}_{i,l}\dot Z_{i,l}=-\bigg(\sum_{j=1}^n
\pi_{i,j}^{l}v'_jZ_j\bigg)W_{i,l}\label{inverse-driver}
\ee
Thus, the aggregate $W$ satisfies another $\sigma\pi$-ODE that we 
call {\it the inverse Driver} (associated to a given original system) \footnote{We point out that the driver equations (\ref{dri}) are quadratic differential equation of the Bernoulli type (Riccati type with no zero degree term). The ODE (\ref{inverse-driver}) is in fact consistent with 
the basic properties of this kind of equations.} 
In accordance with the terminology introduced before, the inverse Driver  is the sub-system $(\cal W,\cal I_{W})$ of the 
ODE (\ref{fin}), (\ref{inverse-driver}) where $\cal I_{W}$ is the set of indices addressing (\ref{inverse-driver}) and $\cal W$ is defined as 
\be
\begin{array}{rl}
{\cal W}&=\{(x,W,t)\in \cal D_r\times \re^d\times {\cal T}\quad /\quad W_{i,l}=Z^{-1}_{i,l},\quad \!\! Z=\Phi(x,t)\},\\
 &={\cal D}_r\times \Phi^W({\cal S}_r)\times \cal T,
\end{array}
\label{SWst}
\ee
where $\Phi^W:{\rm Dom}\{\Phi^W\}={\rm Dom}\{\Phi\}=\cal S_r\subset \re^{n+1}\to \re^d$  has components $\Phi^W_{i,l}$ given by
\be
\Phi^W_{i,l}(x)=(\Phi_{i,l}(x,t))^{-1},\quad x\in \cal D_r.\label{dfPhist}
\ee
We see that the inverse Driver, as well as its 
{\it direct} counterpart, entails the calculation of the original solution, 
provided it feeds an {\it inverse final} system. Indeed, any (regular) solution $(x,t_0)\in \cal D_r\times \cal T$ is given by 
\be
\dot x_i=\Big(\sum_{l=1}^{\nu_i} W_{i,l}^{-1} v_{i,l}\Big)x_i,\qquad
    x_i(t)=e^{\sum_{l=1}^{\nu_i}\int_{t_0}^t W_{i,l}^{-1}v_{i,l}} x_i.\label{inverse-final}
\ee
We name the system $\cal W$ associated to the ODEs (\ref{inverse-driver}), (\ref{inverse-final}), an {\it inverse quadratization} of 
the original $\sigma\pi$-system.\\

\subsubsection{Bernoulli differential equations.} In the second example of \S \ref{sc:SCO} we have met a scalar Driver ODE belonging to a known class of scalar differential equations, 
the {\it Bernoulli differential equations}, that are all solvable through a change of variables that turns them into linear equations. 
Recall that a Bernoulli ODE, in our notation, is a $\sigma\pi$-ODE with $n=1$ (it is scalar), $\nu=2$, $X_{11}=x$, $X_{12}=x^\alpha$ with $1\ne \alpha\in \re$:
\be
\dot x=v_1 x+v_2 x^{\alpha}.\label{BEs}
\ee
Following the classic method, one divides (\ref{BEs}) by $x^{\alpha}$ and obtains
$
x^{-\alpha} \dot x= v_1 x^{1-\alpha} +v_2
$.
Then, by defining $w=x^{1-\alpha}$, one has $\dot w=(1-\alpha) x^{-\alpha} \dot x$, and finally
\be
\dot w=(1-\alpha)(v_1w+v_2).\label{BEw}
\ee
Therefore, $w$ can be found as solution of the linear ODE (\ref{BEw}) and $x=w^{1/(1-\alpha)}$. 

We get the same result (even relaxing the hypothesis $\alpha\ne 1$) by applying to (\ref{BEs}) the exact quadratization in {\it inverse form} given by the inverse Driver (\ref{inverse-driver}). For,  the Driver components associated to the $\sigma\pi$-ODE (\ref{BEs}) are $Z_{11}=1$, and 
$Z_{12}=x^{\alpha-1}$, thus the inverse Driver components are $W_{11}=1$, $W_{12}=x^{1-\alpha}$, and we see that $W_{12}=w$. We get the ODE for $W_{12}$ directly from its monomial, $x^{1-\alpha}$, by applying the inverse-Driver formula (\ref{inverse-driver}) : 
$\dot W_{12}= (1-\alpha)(v'Z)W_{12}=(1-\alpha)(v_1Z_{11}+v_2Z_{12})W_{12}$, then, since $Z_{11}=1$, and $Z_{12}W_{12}=1$ we have 
\be
\dot W_{12}= (1-\alpha)(v_1W_{12}+v_2),\label{IDex}
\ee
which is the linear equation (\ref{BEw}).  For $\alpha=1$ the original ODE is linear, thus solvable with the known classic methods, and unlike Bernoulli method, that fails for $\alpha=1$, such 
solution is still a local sub-solution of the inverse quadratization. 

We can say that a $\sigma\pi$-ODE is a generalized version, in multiple indeterminates, of the Bernoulli equation. 
Unfortunately though, besides the Bernoulli equation, 
there are very few particular cases, and not really meaningful, for which the inverse quadratization 
leads to the solution calculation via linearization.\\ 

 \subsection{The Driver of a linear ODE} 

A linear ODE  is a $\sigma\pi$-ODE where all $p_{i,j}^l$ in (\ref{pi}) may be one or zero, and for any $i,l$ there is one and only one $j$ such that $p_{i,j}^l=1$. For a linear ODE we can always assume that, for any $i$,  
$\nu_i=n$, which amounts to saying that all the $n$ monomials: $X_{i,1}=x_1,\ldots, X_{i,n}=x_n$ of the $i$-th equation 
are actually written (otherwise, write it multiplied by a zero coefficient). This in turn entails that the monomials do not depend of the equation index $i$. In such a framework we have coefficients 
$v_{i,j}$,  for $i,j=1,\ldots,n$ and $p^l_{i,j}=\delta_{j,l}$ for $i,j,l=1,\ldots,n$, (not depending of $i$).  That said, 
consider a linear system: 
\be
\dot x_i=\sum_{j=1}^n v_{i,j}x_j, \quad {\rm \Bigg(} \dot x_i=\sum_{j=1}^{n} v_{i,j}X_{i,j}\quad {\rm with}\quad \!\! X_{i,j}=x_j {\rm \Bigg)}.\label{aff}
\ee
Since (\ref{aff}) is a $\sigma\pi$-ODE, it undergoes a quadratization, and then it has an associated Driver, that we readily compute by the general Driver expression (\ref{dri}).  By (\ref{pi}): $\pi^l_{i,j}=\delta_{j,l}-\delta_{i,j}$, which replaced in (\ref{dri}) yields
\be
\dot Z_{i,l}=(v'_lZ_l-v_i'Z_i)Z_{i,l} 
\label{ldri}
\ee
The Driver lives on the domain ${\cal S}^*\subset \re^{n+d+1}$ given by (\ref{sbdom}) with $\Phi$ given by 
\be
Z_{i,l}=\Phi_{i,l}(x,t)={x_l \over x_i},\quad {\rm for}\quad \!\!\! l=1,\ldots,n. \label{Zaff}
\ee
As explained earlier, we get an inclusive quadratization by enlarging the Driver (\ref{ldri}) with the further $n$ ODEs:
\be
\dot Z_{i,n+1}=(v_i'Z_i)Z_{i,n+1},\label{IDL}
\ee
with $Z_{i,n+1}=x_i$, which replaces the final equation. 

\subsection{Quadratic ODEs, Drivers and self-Drivers}
\label{SD:sec}

We name {\it Driver-type ODE}  (in $n$ indeterminates) an  homogeneous quadratic ODE of the type 
\be
\dot x_i=\sum_{j=1}^{n} v_{i,j}x_ix_{j}\qquad  \Big(=(v'_i x) x_i\Big) \label{SD}
\ee
where the $v_{i,j}$ are real analytic functions all defined on some open interval $\cal T\subset \re$.  Therefore, 
a {\it Driver-type system} shall be the system $\re^n\times \cal T$ associated to the ODE (\ref{SD}). 

This name is motivated by the following theorem. \\

\noindent
\begin{theorem}\label{th:Drty} {\it  Any Driver-type ODE has always one sub-system that is the Driver of some $\sigma\pi$-system
}\end{theorem}

{\it Proof.} It is sufficient to prove that (\ref{SD}) is the Driver ODE of some $\sigma\pi$-ODE. For, let us choose an integer $m\le n$, and $m$ integers $\nu_1,\ldots,\nu_m$ such that $\nu_1+\ldots,\nu_m=n$. Then we can  rewrite (\ref{SD}) as the following $m$-partite system with $k$th part $\nu_k$ dimensional, for $k=1,\ldots,m$:
\be
\dot x_{k,s}=\sum_{r,l=1}^{m,\nu_r} v_{k,s,r,l} x_{k,s}x_{r,l},\label{mpart} 
\ee
where in the summation $r=1,\ldots,m$ runs first, and $s=1,\ldots,\nu_r$, and we have simply replaced $i$ with $(k,s)$, and $j$ with $(r,l)$. 
Now, let us decompose $v_{k,s,r,l}$ as a product of real numbers $\pi^s_{k,r}$ and $v^*_{r,l}$:
\be
v_{k,s,r,l}=\pi^s_{k,r}v^*_{r,l},\label{mpartT1}
\ee
and define the real number $p^s_{k,r}$:
\be
p^s_{k,r}=\pi^s_{k,r}+\delta_{i,j}.
\ee
By using (\ref{dri}), and (\ref{pi}) one readily verify that (\ref{mpart}) ---  with $x$ replaced by the symbol $Z$   --- 
is the Driver ODE associated to the 
$\sigma\pi$-ODE with system function given by (\ref{SP-form1}), where $m$ replaces $n$ and $v^*$ replaces $v$. \hfill$\bullet$\\

Notice, by the proof of the theorem, that there are indeed many different $\sigma\pi$-systems generating  the same Driver ODE 
(\ref{SD}). 
Nonetheless, they generate each a different Driver as a different sub-system of the (same) Driver ODE, in other words, as a set of sub-solution (in the $Z$ indeterminates) lying in a different domain $\cal S^*$ --- which is in fact realized by (\ref{sbdom}) through a different map $\Phi^*$. 

We name {\it self-Driver} (in $n$ indeterminates) any sub-system of the Driver-type ODE (\ref{SD}). Thus, a sub-set $\cal F\subset \re^n\times \cal T$ invariant with respect to the ODE (\ref{SD}), determines a self-Driver, and the ODE (\ref{SD}) itself is a self-Driver (with domain $\re^n\times \cal T$). The set of all self-Drivers associated to the ODE (\ref{SD}) that are Drivers (of some pre-Driver) is a sub-set (in general proper) of the set of all self-Drivers. Clearly, a self-Driver $\cal F\subset \re^n\times \cal T$ is a Driver if and only if there exists an integer 
$m\le n$, $m$ integers $\nu_1,\ldots,\nu_m$ such that $\nu_1+\ldots+\nu_m=n$, and $n$ monomial functions $X_{i,l}$, for $i=1,\ldots,m$ and 
$l=1,\ldots,\nu_i$, such that $\cal F=\Phi(\cal D_r\times \cal T)$ with $\Phi_{i,l}=x_i^{-1}X_{i,l}$. 

Since (\ref{SD}) is quadratic, it is a $\sigma\pi$-ODE as well, of order $n$ and with $n^2$ monomials. 
Thus, it undergoes a canonical quadratization and has in turn an associated Driver. 
The name 'self-Driver' comes from the fact that any system associated with the ODE (\ref{SD}) has a Driver 
whose solutions include, as sub-solutions, all the solutions of the system itself, and only these. 
Indeed, we can prove the following theorem. \\

\noindent
\begin{theorem}\label{th:dristrut} {\it  Let us consider a self-Driver $\cal F\subset \re^n\times \cal T$, associated to the ODE (\ref{SD}).  
Then the Driver of $\cal F$ is defined by the following ODE $n$-partite 
in the $n^2$ indeterminates $Z_{i,j}$, $i,j=1,\ldots,n$:
\be
\dot Z_{i,j}=(v'_jZ_j)Z_{i,j},\label{SSD}
\ee 
on the domain $\cal F^{\bf d}\subset \re^{n^2}\times \cal T$:  
\be
\cal F^{\bf d}=\bigcup_{(x,t)\in {\cal F}}\left\{(Z,t)\in \re^{n^2}\times \cal T\quad /\quad Z_{i,j}=x_i,\quad {\rm for} \quad\!\! 
(x,t)\in \cal F\right\}.\label{DomSD}
\ee
}\end{theorem}

{\it Proof.}  Let us apply (\ref{dri}) to (\ref{SD}).  We have $\pi^l_{i,j}=p^l_{i,j}-\delta_{i,j}$, where $p^l_{i,j}$ is the power of 
$x_j$ in the monomial $x_ix_l$, and hence is $1$ for $j=l$ or $j=i\ne l$. It follows that $\pi^l_{i,j}=1$ if and only if $j=l$. Therefore,   
the Driver ODE is (\ref{SSD}). 
The Driver map $\Phi$, whose general formula is (\ref{Ph}), applied to the $\sigma\pi$-ODE (\ref{SD}) gives: 
$Z_{i,l}=\Phi(x)=x_{l}$, and hence (\ref{DomSD}) is an invariant set for the ODE (\ref{SSD}), thus
the Driver domain.  \hfill$\bullet$\\

Although the Driver of the self-Driver $\cal F$ has $n^2$ indeterminates, 
any solution has many equal components, so that just suitably picking up $n$ among the $n^2$ solution components, 
one get the original self-Driver solution $(x,t)\in \cal F$. This motivates the following identification, stated by the definition below, 
of a self-Driver with its Driver, which amounts to consider the latter just a different (augmented) {\it representation} of the former.\\

\noindent
\begin{definition}\label{df:augR} {\it  Given a self-Driver $\cal F\subset \re^n\times \cal T$ described by the ODE (\ref{SD}), we call 
the system (in $\re^{n^2}\times \cal T$)  
defined by the ODE (\ref{SSD}) on the domain $\cal F^{\bf d}$ given by (\ref{DomSD}) the  
'augmented representation' of the self-driver $\cal F$.}
\end{definition}
\vskip0.8cm

We finally give the following definition, acconding to which, s Driver-type ODE is characterized by its own matrix of coefficients $v_{i,j}$-\\  

\noindent
\begin{definition}\label{df:frame} {\it We name 'the frame' of a Driver-type ODE in $n$ indeterminates,  
the matrix $V\in\re^{n\times n}$, $V=(v_{i,j})$, collecting all of its own coefficients. 
}
\end{definition}
\vskip0.8cm

\section{The solution formula for Driver-type ODEs}
\label{sc:SDTheo}

\subsection{Statement of the main theorem}
\label{ssc:mai}

The importance of Driver type quadratic ODEs, and, hence, of the exact quadratization, 
relies on the result stated in the next theorem, showing that  
all solutions of a Driver type ODE 
can be calculated through power series. 

Let us consider an ODE of the type (\ref{SD}) with analytic coefficients defined on the interval $\cal T\subset \re$. The associated 
{\it full} system is then $(\re^n\times \cal T, \{1,\ldots,n\})$, but on account of the properties of analytic systems, described in 
\S \ref{ssc:sinreg}, without loss of generality we can consider the sub-system 
$(\cal D\times \cal T, \{1,\ldots,n\})$, where $\cal D=\re^n\setminus S^0$, and $S^0$ is the set of all coordinate  hyperplanes in $\re^n$. 
Indeed, any solution $(x,t)\in \re^n\times \cal T$ with $x_i=0$ for some $i$, is either a regular solution, and, hence, it can be rewritten as 
$(x^*,t^*)\in \cal D\times \cal T$, or is singular, and then it can be rewritten as $(x^*,t)\in \cal D^*\times \cal T$, where $\cal D^*\times \cal T$ 
is the domain of system projected on $\{1,\ldots, n\}\setminus \{i\}$. Besides, if the system is a self-Driver, its domain, say $\cal D\times \cal T\subset \re^n\times \cal T$, is the image: $\cal D=\Phi(\cal D^*)$ where  $\cal D^*\in \re^m$, for some $m\le n$,  
$\cal D^*\cap D^0=\emptyset$ ($D^0$ denoting now the set of points of all hypherplanes in $\re^m$), and $\Phi$ is an $\re^n$-valued 
map whose components are all {\it monomials} of the type $z^{p_1}_1\dots z^{p_m}_m$ with $z_1,\ldots, z_m$ real indeterminates, and 
$p_1,\ldots p_m$ real exponents. Thus, any component of $\Phi$ can be zero only on $D^0$, which implies that $\cal D\cap S^0=\emptyset$. 
The hypothesis of a domain not including any point of the coordinate hypherplanes, is then directly verified for the more expected situation 
of our study: i.e. when the ODE is a quadratization of some original $\sigma\pi$-reducible system.\\

With that being said, for the solution $(x,t_0)$ of the ODE (\ref{SD}) with coefficients $v_{i,j}$, 
where $(x,t_0)\in \cal D\times \cal T$, we can consider the Taylor expansion of the corresponding 
curve-solution, $x(t)$, with $t\in {\bf I}_{x,t_0}\subset \cal T$: 
\be
x_i(t)=\sum_{k=0}^{\infty}c_k(i) {(t-t_0)^k\over k!},\quad c_k(i)=x^{(k)}(t_0), \label{SPSD}
\ee
where $x^{(k)}(\cdot)$, defined for all $t\in {\bf I}_{x,t_0}$, denotes the curve of the $k$th time-derivative of $x(\cdot)$. 
Identity (\ref{SPSD}) holds for any $t\in {\bf I}_{x,t_0}$ such that the series in the right hand side converges, which happens in general 
on a sub-interval:
\be
{\bf I}(t_0,r)\subset {\bf I}_{x,t_0},\label{intconv}
\ee
where ${\bf I}(t_0,r)$ is the open interval centered in $t_0$ with ray $r>0$. 

The following theorem gives the general formula for the coefficients $c_k(i)$. In the following the partial time-derivative operator 
$\partial_t$ will be used in place of the 'dot' notation for symbols including superscripts, thus $\partial_t v^k$ in place of $\dot {\bf v}$ for 
${\bf v}=v^k$.\\ 

\noindent
\begin{theorem}\label{th:serie} {\it The coefficients $c_k(i)$ in the Taylor expansion (\ref{SPSD}) are given by 
\bea
&&c_0(i)=x_i,\label{GE0}\\
&&c_k(i)=\sum_{s=2}^{k+1}\sum_{i_1,\ldots,i_{s-1}\in \cal I} v^{k+1,s}_{i,i_1,\ldots,i_{s-1}}(t_0) \cdot x_ix_{i_1}\cdots x_{i_{s-1}}, \quad \cal I=\{1,\ldots,n\},\quad \forall k\ge 1,\label{GEs}
\eea
where $i_l\in \cal I$, for $l=0,1,\ldots, k$ (we set $i_0=i$), and the coefficients $v^{k+1,s}_{i,i_1,\ldots,i_{s-1}}$, 
($=v^{k+1,s}_{i,i_1,\ldots,i_{s-1}}(t)$,  $t\in {\bf I}_{x,t_0}$) are given 
by the following recursive equations ($k\ge 2$):
\bea
&&v^{k+1,k+1}_{i,i_1,\ldots,i_k}=v^{k,k}_{i,i_1,\ldots,i_{k-1}}\left(\sum_{j=0}^{k-1} v_{i_j, i_k}\right), \label{Vkk}\\
&&v^{k+1,s}_{i,i_1,\ldots,i_{s-1}}=v^{k,s-1}_{i,i_1,\ldots,i_{s-2}}\left(\sum_{j=0}^{s-2} v_{i_j, i_{s-1}}\right)+\partial_t v^{k,s}_{i,i_1,\ldots,i_{s-1}},\quad s=3,\ldots,k,\label{Vks}\\
&&v^{k+1,2}_{i,i_1}=\partial_t v^{k,2}_{i,i_1}, \hskip 2cm \Big( = \partial_t^{k+1} v_{i,i_1}\Big),\label{Vk2}
\eea
which are initialized by 
\be
v^{2,2}_{i,i_1}=v_{i,i_1}.\label{initM}
\ee
}
\end{theorem}
\vskip0.8cm

\subsubsection{Comments to the statement of Theorem \ref{th:serie}}
\label{ssc:com} for any $k,s$, the number of coefficients $v^{k,s}_{i,i_1,\ldots,i_{s-1}}$ is $n^s$. in order to calculate 
$c_k(i)$, $i$ is fixed, and one needs $n+\ldots +n^k$ coefficients, one for any fixed sequence 
$i_1,\ldots,i_{s-1}\in \cal I$, $s=2,\ldots, k+1$. The coefficients have to be multiplied each for the corresponding monomial $x_ix_{i_1}\cdots x_{i_{s-1}}$ obtained by the components of 
the initial value $x=x(t_0)$. By using the index convention we have $v^{k,s}_i(t)\in \re^{n^s-1}$, and we can write (\ref{GEs})
 in vector notation as well:
$$
c_k(i)=\sum_{s=2}^{k+1}\left(v_i^{k,s} (t_0)\right)' x^{\otimes s-1} x_i
$$
where $x^{\otimes l}$ denotes the $l$-times tensor product $x\otimes \ldots\otimes x$. 
Note that the scheme (\ref{Vkk})--(\ref{initM}) is recursive in $k$: from the set of all coefficients involved 
in $c_{k-1}(i)$, namely: 
$$
\cal V_i(k-1)=\{ v^{r,s}_i\quad |\quad r=2,\ldots,k,\quad s=2,\ldots, r\}
$$
one directly calculates from (\ref{Vkk})--(\ref{Vk2}) the set $\{v_i^{k+1,s}, \quad | s=2,\ldots, k+1\}$, which,  joined to $\cal V_i(k-1)$, 
yields $\cal V_i(k)$. 
Note that the coefficients $v^{k,s}_i$ can be calculated independently of the initial point $x$. 
If formula (\ref{SPSD}) is used numerically, then for any fixed approximation degree, and initial point $x\in \cal D$, it allows to calculate separately each component $x_i(t)$ of the solution, within the given approximation, 
and within the convergence interval ${\bf I}(t_0,r)$. One is able as well to calculate a possible extension of the solution outside ${\bf I}(t_0,r)$, by moving the initial time to another $t_0'\in {\bf I}(t_0,r)$ -- provided that $x(t'_0)\in \cal D$, and  ${\bf I}(t_0',r')\cup {\bf I}(t_0,r)\setminus {\bf I}(t_0,r)\ne \emptyset$ -- and calculating the value at $t_0'$ of always the same unique pre-calculated set of coefficients $v_i^{k,s}$.

\subsection{A preliminary result}
\label{ssc:prel}

In order to prove Theorem \ref{th:serie} we need some additional definitions. Let us fix $i\in \cal I$, 
where $\cal I$ is the set of indices in (\ref{GEs}). 
For any integer $k\ge 2$, we define $\cal P_i^k$ to be the set of all homogeneous polynomials of degree $k$ of the $n$ indeterminates 
$x_1,\ldots,x_n$, with $x_i$ always in, and where the coefficients are  time-functions, defined on some open real interval, and having time-derivative of any order. Therefore, $p_i^k\in \cal P_i^k$ has the form:
\be
p_i^k(x)=\sum_{i_1,\ldots,i_{k-1}\in \cal I} v^k_{i,i_1,\ldots,i_{k-1}} x_ix_{i_1}\ldots, x_{i_{k-1}}.\label{pk}
\ee
Now, let us consider a Driver-type ODE with frame $V$, and, related to $V$, define a map $A_i:\cal P_i^k\to \cal P_i^{k+1}$, for any $k\ge 2$ as follows:
\be
A_ip_i^k=p_i^{k+1}=\sum_{i_1,\ldots,i_k\in \cal I} v^{k+1}_{i,i_1,\ldots,i_{k}} x_ix_{i_1}\ldots, x_{i_{k}},\quad {\rm with}\quad \!\! 
v^{k+1}_{i,i_1,\ldots,i_{k}}=v^{k}_{i,i_1,\ldots,i_{k-1}}\left(\sum_{j=0}^{k-1}v_{i_j,i_k}\right),\label{dA}
\ee
where $v_{i_j,k}$ are coefficients of the frame $V$, and conventionally we set $i_0=i$. Note that any $\cal P_i^k$ is a vector space, with naturally defined operations and null vector, thus we can define $\cal P_i$ as the (infinite) direct sum of all $\cal P_i^k$'s:
$
\cal P_i=\oplus_{s=0}^{+\infty} \cal P_i^s$. 
Therefore, for any $q\in\cal P$ there is a $k\ge 2$ and $p_i^s\in \cal P_i^s$ for $s=2,\ldots,k$ such that $q_i=\sum_{s=2}^kp^s$, and we say that $k$ is the degree of $q_i$. We can naturally extend as follows the definition of $A_i$ given earlier, so that $A_i:\cal P\to \cal P$: 
let $q_i^k\in \cal P_i$ a $k$-degree polynomial, then $q_i^k=p_i^2+\ldots+p_i^k$, for some $p_i^s\in \cal P_i^s$, then we define
$
A_iq_i^k=\sum_{s=2}^k A_ip_i^s
$. The exponent $A_i^s$ is defined at once as a composition of maps, and we have $A_i^s: \cal P_i^k\to \cal P_i^{k+s}$ for any integer $s\ge 0$, provided  we set $A_i^0=I$ (the identity map). Besides $A_i$ we make use of the partial time-derivative operator $\partial_t$, for which a similar 
convention is applied: $\partial_t^s$, for $s\ge 0$ shall denote composition, and $\partial^0_t$ is assumed to be the identity. Note that 
$\partial^s_t: \cal P_i^k\to \cal P_i^k$ for any $k\ge 2$, $i\in \cal I$, and we have:
\be
\partial^s_tp_i^k=\sum_{i_1,\ldots,i_{k-1}} \partial_t v^{k}_{i,i_1,\ldots,i_{k-1}} x_ix_{i_1}\ldots, x_{i_{k-1}}.\label{papk}
\ee
In the same way all compositions of the maps $A$ and $\partial_t$ are well defined, as well as the map $(A+\partial_t)^s$ for any $s\ge 0$, and 
we have $(A_i+\partial_t)^s \cal P_i\subset \cal P_i$, and $(A_i+\partial_t)^s \cal P_i^k\subset \cal P_i^{k+s}$. 

Now, hereinafter we fix the symbol $p^2_i\in \cal P^2$ to be the polynomial whose coefficients are just the components of the frame $V$:
\be
p^2_i(x)=\sum_{i_1} v_{i,i_1} x_i(t)x_{i_1}(t), \label{p2}
\ee
which, in the formalism defined above means that we are taking $v^{2}=v$. It follows at once that we can rewrite a Driver-type ODE with frame $V$ as:
\be
x_i^{(1)}=p^2_i(x). \label{DTOrwt}
\ee
The proof of Theorem \ref{th:serie} requires the following result, that we single out in a lemma because it has  
a significance in its own right. \\ 

\noindent
\begin{lemma}\label{lm:serie} {\it There exists an open interval 
\be
{\bf I}^*_{x,t_0}\subset {\bf I}_{x,t_0},\label{opin}
\ee
with $t_0\in {\bf I}^*_{x,t_0}$, such that 
the $k$th time-derivative of the curve $x(\cdot)$ of the solution $(x,t_0)\in \cal D\times \cal T$ of (\ref{SD}), is given component-wise by 
\be
x_i^{(k)}(t)=q_i^{k+1}(x(t))\hskip 1cm \forall k\ge 1,\qquad \forall t\in {\bf I}^*_{x,t_0},\label{kder}
\ee
where  $q_i^{k}\in \cal P_i$, is a $k$th degree polynomial for any $k\ge 2$ satisfying the recursive equation:
\be
q_i^{k+1}=(A_i+\partial_t) q_i^k, \qquad q_i^2=p_i^2.\label{qkrec}
\ee
}
\end{lemma}
\vskip0.8cm

{\it Proof.} Since $t\mapsto x(t)$ is continuous in ${\bf I}_{x,t_0}$, and $x_i\ne 0$ $\forall i$, there exists an interval 
${\bf I}^*_{x,t_0}\subset {\bf I}_{x,t_0}$, with $t_0\in {\bf I}^*_{x,t_0}$, such that 
\be
x_i(t)\ne 0, \quad \forall t\in {\bf I}^*_{x,t_0}.\label{xine0}
\ee
Note that the recursive equation (\ref{qkrec}) is equivalent to
\be
q_i^{k+1}=(A_i+\partial_t)^{k-1}p_i^2,\label{qnonrec}
\ee
therefore, (\ref{kder}) is equivalent to
\be
x_i^{(k)}(t)=(A_i+\partial_t)^{k-1}p_i^2(x(t)), \hskip 1cm \forall k\ge 1,\qquad \forall t\in {\bf I}^*_{x,t_0}.\label{kderts}
\ee
We prove the lemma by induction. 
For $k=1$ the lemma holds, since (\ref{kderts}) turns into (\ref{DTOrwt}).  
Suppose (\ref{kderts}) is true, and let us prove (\ref{kderts}) 
for $k+1$ replacing $k$. We have  ($D_t$ denoting the total time-derivative, replacing the {\it dot} notation upon symbols with 
other superscripts):
\be
x_i^{(k+1)}(t)=D_{t} x_i^{(k)} (t)=\sum_{l=1}^n \left(\partial_{x_l} x_i^{(k)}\right)\Big|_{x=x(t)} \dot x_l(t)+ \partial_t x^{(k)}_i
\Big|_{x=x(t)}\label{lmM1}
\ee
As for the last term, by using the inductive hypothesis, it results in
\be
\partial_t x^{(k)}\Big|_{x=x(t)}=\partial_t(A_i+\partial_t)^{k-1} p^2_i(x(t)).\label{lmM2}
\ee
By  (\ref{qnonrec}), $q_i^{k+1}$ is a polynomial of degree $k+1$ in $\cal P_i$, therefore it has the form:
\be
q_i^{k+1}=\sum_{s=2}^{k+1}p^s_i, \qquad p^s_i\in \cal P^s,\label{dq}
\ee
where the $p_i^s$'s are homogeneous polynomials of the type (\ref{pk}).  
Thus, the first term in the left hand side of (\ref{lmM1}) becomes, by replacing the polynomial expression (\ref{dq}) for $q^{k+1}_i$:
\be
\sum_{l=1}^n \left(\partial_{x_l} x_i^{(k)}\right)\Big|_{x=x(t)} \dot x_l(t)
= \sum_{s=2}^{k+1}\sum_{i_1,\ldots,i_{s-1}} \left(v^{s}_{i,i_1,\ldots,i_{s-1}} \sum_{l=1}^n \partial_{x_{l}}  (x_ix_{i_1}\cdots x_{i_{s-1}})\Big|_{x=x(t)}\cdot x_i(t)\right)
\label{dxtmp}
\ee
Now, on account of (\ref{xine0}), we have
\bea
\sum_{l=1}^n \partial_{x_{l}}  (x_ix_{i_1}\cdots x_{i_{s-1}})\Big|_{x=x(t)}
\cdot \dot x_l(t) &&= \sum_{j=0}^{s-1} \partial_{x_{i_j}}  (x_ix_{i_1}\cdots x_{i_{s-1}}) \Big|_{x=x(t)}\cdot \dot x_{i_j}(t), \nonumber\\
&&= \sum_{j=0}^{s-1}   x_ix_{i_1}\cdots x_{i_k} x^{-1}_{i_j} \Big|_{x=x(t)}
\cdot \dot x_{i_j}(t),\nonumber
\eea
thus, by replacing $\dot x_{i_j}$ with the right hand side of (\ref{p2}) where we also rename $i_1$ after the new index $i_{s}$, we can write
\bea
&&\hskip-1.5cm \sum_{i_1,\ldots,i_{s-1}} \left(v^{s}_{i,i_1,\ldots,i_{s-1}} \sum_{j=0}^{s-1}   x_ix_{i_1}\cdots x_{i_{s-1}} x^{-1}_{i_j}
\cdot 
\sum_{i_{s}} v_{i_j,i_{s}} x_{i_j}x_{i_{s}},\right)\Bigg|_{x=x(t)} \nonumber\\
&&=\sum_{i_1,\ldots,i_{s}} \left(v^{s}_{i,i_1,\ldots,i_{s-1}}   x_ix_{i_1}\cdots x_{i_{s-1}}  
\cdot 
\sum_{j=0}^{s-1} x^{-1}_{i_j} v_{i_j,i_{s}} x_{i_j}x_{i_{s}},\right)\Bigg|_{x=x(t)}\nonumber\\
&&=\sum_{i_1,\ldots,i_{s}} \left(v^{s}_{i,i_1,\ldots,i_{s-1}}   x_ix_{i_1}\cdots x_{i_{s}}\right) \Big|_{x=x(t)}\cdot 
\left(\sum_{j=0}^{s-1} v_{i_j,i_{s}}\right).\nonumber\\
&&=A_ip^{s}(x(t)), \nonumber
\eea
where, in the last step, the definition of the map $A_i$, given in (\ref{dA}) has been used. 
By replacing the above identity in (\ref{dxtmp}), and using again the inductive hypothesis -- in the form (\ref{kder}) -- 
we obtain
$$
\sum_{l=1}^n \left(\partial_{x_l} x_i^{(k)}\right)\Big|_{x=x(t)} \dot x_l(t)= \sum_{s=2}^{k+1}A_ip^s_i=A_iq^{k+1}_i=
A_i(A_i+\partial_t)^{k-1}p^2_i(x(t))
$$
which used in (\ref{lmM1}) at once with (\ref{lmM2}) gives 
$$
x_i^{(k+1)}(t)=(A_i+\partial_t)^k p^2_i(x(t)),
$$
Q.E.D.\hfill$\bullet$\\ 

\subsection{Proof of Theorem \ref{th:serie}.} 
\label{scc:pro}
 By Lemma \ref{lm:serie} the coefficients $c_k(i)$ are given by 
\be
c_k(i)=x^{(k)}(t_0)=q^{k+1}_i(x), \qquad k\ge 0\label{ck}
\ee
and since $q^{k+1}_i$ is a $k+1$ degree polynomial in $\cal P_i$,  
there exist $k-1$ polynomials $p^{k+1,s}_i\in \cal P_i^{s}$, $s=2,\ldots,k+1$ such that 
\be
q^{k+1}_i=\sum_{s=2}^{k+1}p^{k+1,s}_i=\sum_{s=2}^{k+1} \sum_{i_1,\ldots,i_{s-1}} v^{k+1,s}_{i,i_1,\ldots,i_{s-1}} 
x_ix_{i_1}\ldots x_{i_{s-1}}\label{genfpk}
\ee
which gives the general expression (\ref{GEs}). In order to proof formulas (\ref{Vkk})--(\ref{Vk2}), consider the recursive formula 
(\ref{qkrec}). where $q^k_i$ has a similar expression like (\ref{genfpk}) but the degree is $k$: 
\be
q^{k}_i=\sum_{s=2}^{k}p^{k,s}_i=\sum_{s=2}^{k} \sum_{i_1,\ldots,i_{s-1}} v^{k,s}_{i,i_1,\ldots,i_{s-1}} 
x_ix_{i_1}\ldots x_{i_{s-1}}\label{genfpk2}
\ee
thus we have (omit all subscripts)
\be
\begin{array}{rl}
q^{k+1}_i=p^{k+1,k+1}+p^{k+1,k}+\ldots+p^{k+1,2}&=Ap^{k,k}+Ap^{k,k-1}+\ldots,Ap^{k,2}\\
&+\partial p^{k,k}+\partial p^{k,k-1}+\ldots,\partial p^{k,2}                                                                                 
\end{array}\label{q=p}
\ee
where in the right hand side, by definition given in (\ref{dA}), we have that $Ap^{k,s}\in \cal P^{s+1}$, and 
$\partial p^{k,s}\in \cal P^s$, for $s=2,\ldots,k$. By equating the homogeneous polynomials with the same degree, (\ref{q=p}) is equivalent to  
the following set of  equalities:
\bea
p^{k+1,k+1}&=&Ap^{k,k},\nonumber\\
p^{k+1,s}&=&Ap^{k,s-1}+\partial p^{k,s},\qquad s=3,\ldots, k\nonumber\\
p^{k+1,2}&=&\partial p^{k,2}\nonumber
\eea
which directly gives equalities (\ref{Vkk})--(\ref{Vk2}), as soon as the definition of A, given in (\ref{dA}) is applied to the polynomial expression of $p^{k,s}$ and $p^{k+1,s}$, for $s=2,\ldots k+1$, which can be read out from (\ref{genfpk2}) and (\ref{genfpk}).
\hfill $\bullet$\\

Let $\cal S\subset \{1,\ldots, n\}$ be the set of indices: 
\be
\cal S =\{ j\in \{1,\ldots,n\}\quad | \quad  \exists i, \quad \!\!{\rm such} \!\! \quad \!\!{\rm that} \!\! \quad v_{i,j}\ne 0\},\label{defSc}
\ee
In other words, the set $\cal S$ is the set of indices addressing the non-zero columns of the frame $V$. The following theorem can ease the calculation of formula (\ref{GEs}).\\

\noindent
\begin{theorem}\label{th:S} {\it In formula (\ref{GEs}) the indices $i,i_1,\ldots,i_{s-1}$, don't have to run across all values $\{1,\ldots, n\}$: it's enough they span just the subset $\cal S$ defined in (\ref{defSc}).
}
\end{theorem}
\vskip0.8cm

{\it Proof.} The thesis is equivalent to say that for any $k\ge 2$ and $2\le s\le k$, the following condition is satisfied:\\

{\it Condition (C)}:
$v^{k,s}_{i,i_1,\ldots,i_{s-1}}=0$ if  there is an $m$, $1\le m\le s-1$ such that $i_{m}\not \in \cal S$.\\

Let is denote by ${\bf v}(l)$, with $l\in \na$, the set of all pair $(k,s)$ with $k\ge 2$, $2\le s\le k$, and $k-s=l$, then the theorem is proven 
as soon as we prove, by induction over $l\in \na$, that condition (C) is verified by all pairs in $\cup_{l\in \na}{\bf v}(l) $. For, let us 
consider first the pairs of ${\bf v(0)}$, which address the coefficients 
$v^{k,k}$ generated by (\ref{Vkk}). Suppose that 
there is an index, say $i_m$, $1\le m\le s-1< k$, such that 
$i_m\not\in \cal S$.  Then $v_{j,i_m}=0$ (for any $j$) and by formula (\ref{Vkk}) for $k=m$, 
we have 
$$
v^{m+1,m+1}_{i,i_1,\ldots,i_m}=v^{m,m}_{i,i_1,\ldots,i_{m-1}} \left(\sum_{j=0}^{m-1} v_{i_j, i_m}\right)=0,
$$
and, hence, again by the recursive formula (\ref{Vkk}):
$$
v^{k+1,k+1}_{i,i_1,\ldots,i_k}=v^{m+1,m+1}_{i,i_1,\ldots,i_{m}} \left(\sum_{j=0}^{m} v_{i_j, i_{m+1}}\right)\cdots 
 \left(\sum_{j=0}^{k} v_{i_j, i_{k+1}}\right)=0.
$$
Now, suppose that condition (C) is verified by all pairs in ${\bf v}(l)$, and let us prove that it is verified by all pairs of ${\bf v}(l+1)$ as well. 
For, suppose that $k-s=l$ and there is an $i_m$, $1\le m\le s-1< k$, such that 
$i_m\not\in \cal S$, then, by the inductive hypothesis, we have $v^{k,s}_{i,i_1,\ldots,i_{s-1}}=0$, which simplifies the recursive formula (\ref{Vks}) as follows:  
$$
v^{k+1,s}_{i,i_1,\ldots,i_{s-1}}=v^{k,s-1}_{i,i_1,\ldots,i_{s-2}}\left(\sum_{j=0}^{s-2} v_{i_j, i_{s-1}}\right), 
$$
from which, by successive substitutions up to $m+1$ we obtain
$$
v^{k+1,s}_{i,i_1,\ldots,i_{s-1}}=v^{m+l+1,m}_{i,i_1,\ldots,i_{m-1}} \left(\sum_{j=0}^{s-2} v_{i_j, i_{s-1}}\right)\cdots 
 \left(\sum_{j=0}^{m-1} v_{i_j, i_{m}}\right)=0.
$$
since, again by the inductive hypothesis, $v_{j,i_m}=0$ for all $j$. We have proven that $v^{k+1,s}=0$ for $k-s=l$, which means 
$v^{k,s}=0$ for $k-s=l+1$, i.e. for all pairs in ${\bf v}(l+1)$, Q.E.D. \hfill$\bullet$\\

\subsection{The stationary case}
\label{ssc:statcase}

Worth of mention is the {\it stationary} case, i.e. when the coefficients $v_{i,l}$ are all constants, because in this case the calculation of the 
coefficients $c_k(i)$ greatly simplifies, as shown in the following theorem\\

\noindent
\begin{theorem}\label{th:seriest} {\it In the stationary case, the coefficients $c_k(i)$ in the Taylor expansion (\ref{SPSD}) are constant and are given by 
\be
c_k(i)=\sum_{i_1,\ldots,i_k\in \cal S} v^{k+1}_{i,i_1,\ldots,i_k}\cdot x_ix_{i_1}\cdots x_{i_k},\label{cstat}
\ee
where $i_s=1,\ldots, n$, for $s=0,1,\ldots, k$ (we set $i_0=i$), and the (constant) coefficients $v^{k+1}_{i,i_1,\ldots,i_k}$, 
are given 
by the following recursive equation 
\be
v^{k+1}_{i,i_1,\ldots,i_k}=v^{k}_{i,i_1,\ldots,i_{k-1}}\left(\sum_{j=0}^{k-1} v_{i_j, i_k}\right), \quad v^1_i=1.\label{Vkstat}
\ee
}
\end{theorem}
\vskip0.8cm

{\it Proof.} Let us prove that $v^{k,s}=0$, for $k-s>0$ (note that $k-s<0$ doesn't apply). The proof is by induction over $k-s$, thus consider first the case $k-s=1$, which 
amounts to prove that $v^{k+1,k}=0$, $\forall k>2$. For, let us proceed by induction over $k$: for $k=2$ by (\ref{Vk2}) we have 
$v^{3,2}=\partial_t v^{2,2}=0$ by the stationarity hypothesis. Suppose $v^{k,k-1}=0$, then (\ref{Vks}) for $s=k$ yields 
$v^{k+1,k}=\partial _t v^{k,k}=0$, since by (\ref{Vkk}) initialized by (\ref{Vk2}), and by stationarity, we have $\partial_t v^{k,k}=0$. 
Coming back to the main induction, let us suppose that for $k-s=l$, $v^{k,s}=0$. Then by (\ref{Vks}): 
\be
v^{k+1,s}_{i,i_1,\ldots,i_{s-1}}=v^{k,s-1}_{i,i_1,\ldots,i_{s-1}}\left(\sum_{j=0}^{s-2}v_{i_j,i_{s-1}}\right),\label{k-s=l}
\ee
using which, we can prove that $v^{k,s}=0$ for $k-s=l+1$, by induction over $s$: for $s=2$ we have $k=l+3$, and $v^{l+3,2}=0$ by 
(\ref{Vk2}) and by stationarity. Then, suppose $v^{k,s-1}=0$, and $k-s+1=l+1$, then $k-s=l$ and by (\ref{k-s=l}) we have 
$v^{k+1,s}=0$ and $k+1-s=l+1$, which completes the induction for $k-s>0$. In conclusion, there remains only the coefficients 
$v^{k,k}$ given by (\ref{Vkk}), initialized by (\ref{initM}): the thesis follows, by renaming $v^{k,k}\to v^k$, and extending $v$ with  
$v^1_i=1$ -- which makes (\ref{Vkstat}) returning (\ref{initM}) for $k=1$. 
\hfill$\bullet$\\

Let $\alpha^l_{i,i_1,\ldots,i_{m}}$ for $m\in \na$ be the number of times that $l=i_j$ for $j=0,1,\ldots, m$. Thus, given  
the sequence of indices $i,i_1,\ldots,i_{m}$, we have recursively:
\be
\alpha^l_{i,i_1,\ldots,i_m} = \alpha^l_{i,i_1,\ldots,i_{m-1}}+\delta_{l,i_m}\quad \alpha^l_i=\delta_{l,i},\label{alpharec} 
\ee
and
\be
\alpha^l_{i,i_1,\ldots,i_m}=\sum_{j=0}^m \delta_{l,i_j}\label{alpha}
\ee
It is possible to reduce the computational burden of  formulas (\ref{Vkstat}) and  (\ref{cstat}) by avoiding the calculation of many 
zero terms, as shown in the following corollary. \\

\noindent
\begin{corollary}\label{co:seriest} {\it Formula (\ref{Vkstat}) can be rewritten as follows: 
\bea
&&v^{k+1}_{i,i_1,\ldots,i_k}=v^{k}_{i,i_1,\ldots,i_{k-1}}  \gamma_{i,i_1,\ldots,i_k}  \qquad v^1_i=1,\label{Vkd}\\
&&\gamma_{i,i_1,\ldots,i_k}=\sum_{l\in \rho(i_k)\cup \{i\}}
  \alpha^{l}_{i,i_1,\ldots,i_{k-1}}  
v_{l, i_k}.\label{Vkb}
\eea
where, for any $j\in \cal S$: 
\be
\rho(j) =\{ l\in \cal S \quad | \quad  v_{l,j}\ne 0\}\label{rho}
\ee
}
\end{corollary}
\vskip0.8cm

{\it Proof.} Define 
\be
\tilde \rho(j) =\{l\in  \{1,\ldots,n\}\quad | \quad  v_{l,j}\ne 0\}.\nonumber
\ee
We see in (\ref{Vkstat})  that the index $i_j$ for $j\ne 0$ (thus, with the exception of $i$) doesn't have to take values that are not in 
$\tilde \rho(i_k)$. Since $i_j$, for $j\ne 0$, doesn't have as well to take values that are not in $\cal S$, the thesis follows. \hfill$\bullet$\\ 

Denote by $\sigma$ the cardinality of $\cal S$, $v_M=\max_{i,j\in {\cal I}} |v_{i,j}|$, and for $x\in \cal D$, (thus, $\forall i$, $x_i\ne 0$) 
$x_M=\max_i |x_i|$. In the following Theorem,  we find a lower bound, $\bar r$,  for $r$, the convergence ray of the series for all the functions 
$t\mapsto |x_i(t)|$, components modules of the solution $(x, t_0)$, as well as an upper bound for all of them, in the interval 
$(t_0-\bar r, t_0+\bar r)\subset {\bf I}(t_0,r)$. \\

\noindent
\begin{theorem}\label{th:seriest} {\it 
For all $t\in (t_0-\bar r, t_0+\bar r)\subset {\bf I}(t_0,r)$, with 
\be
\bar r={1\over \sigma v_Mx_M},\label{LB}
\ee
and for any $i=1,\ldots,n$ we have:
\be
|x(t)|\le {x_M\over 1-\sigma v_Mx_Mt}, \label{UB}
\ee
}
\end{theorem}
\vskip0.8cm

{\it Proof.} By (\ref{Vkstat}) we have $|v^{k+1}|\le |v^k| k v_{M}$
which, iterated, yields (on account that $v^1=1$):
$$
|x(t)| \le k! v_M^k.
$$
Thus, by (\ref{cstat}), we have 
$$
|c_k(i)| \le \sigma^k k! v_M^kx_M^{k+1}
$$
which gives a geometric series whose generator is $\sigma v_Mx_M t$, and we get the upper-bound (\ref{UB}). The lower-bound $\bar r$ in (\ref{LB}), is readily derived by an application of the ratio criterion. 
\hfill$\bullet$\\

\subsection{Examples}
\label{ssc:exfin}

\subsubsection{The basic linear ODE} 
\label{sss:blo}
$\dot x=ax$, with $a\in \re$. The solution $(x,0)\in \re^2$ is $x(t)=x e^{at}$.  
Then, let us verify that the Taylor expansion of the solution, calculated by using Theorem \ref{th:seriest} is: 
\be
x(t)=x\sum_{k=0}^\infty {a^k t^k\over k!}.\label{ex-eat}
\ee
For, we have first to quadratize the ODE, which is readily done in this simple case 
by defining the two indeterminates $x_1, x_2$ and setting 
\be
\begin{array}{l}
\dot x_1=ax_1x_2,\\
\dot x_2=0,
\end{array}\label{ex-newe}
\ee
from which we can recover the original solution as the sub-solution (in $x_1$) of the solution $( (x_1,1), 0)$ of (\ref{ex-newe}). Now, the frame of \ref{ex-newe} is
\be
V=\left( 
\begin{array}{cc}
0& a\\
0&0
\end{array}
\right)\label{Vex}
\ee
thus $\cal S=\{2\}$, and by setting $i=1$, the indices $i,i_1,\ldots,i_k$ in (\ref{cstat}) can take just the single value $1,2_{k}$, 
where we denote $2_k=2,\ldots,2$ ($k$-times). Therefore, formula (\ref{Vkstat}) yields:
$$
v^{k+1}_{1,2_k}=v^k_{1,2_{k-1}} (v_{1,2}+kv_{2,2})=v^k_{1,2_{k-1}}a=a^k
$$
and by (\ref{cstat}), recalling that $x_2=1$:
$$
x_1(t)=\sum_{k=0}^\infty v^{k+1}_{1,2_k} x_1 x_2\cdots x_2 {t^k\over k!}=x_1\sum_{k=0}^\infty {a^kt^k\over k!}.
$$
which is (\ref{ex-eat}). \\

\subsubsection{A simple linear non-stationary ODE.} 
\label{sss:blons}Consider the ODE: $\dot x=2t x$, for which the solution $(x,0)$ can be easily calculated as
\be
x(t)=xe^{\int_0^t 2\tau d\tau}=xe^{t^2},\label{ex2-eat2}
\ee
whose Taylor expansion centered at $0$ is:
\be
x(t)=x\sum_{k=0}^\infty {t^{2k}\over k!}=x\left(1+ t^2+{t^4\over 2} + {t^6\over 3!}+\ldots \right),\label{TE}
\ee
To illustrate the point, we now apply Theorem \ref{th:serie} in order to calculate the first four coefficients of the power series in (\ref{TE}).  
Note that we get the same quadratization (\ref{ex-newe}), and frame (\ref{Vex}), as in the previous example, provided $2t$ replaces $a$, 
therefore we have to calculate the sub--solution $x_1$ of solution $((x_1,1),0)$ of the quadratization.  
We have again ${\cal S}=\{2\}$, and $i=1$, thus, on account of Theorem \ref{th:S}, in (\ref{GEs}) the summation in $1,i_1,\ldots, i_{s-1}$ has only one element,  
addressed by the string $1,2_{s-1}$, but the sum in $s$ has $k-1$ elements:
\be
c_k(1)=(v_{1,2}^{k+1,2}+v_{1,2,2}^{k+1,3}+\ldots+v_{1,2_{k}}^{k+1,k+1}) \big |_{t=0}x_1.\label{ex2-gf}
\ee
For $k=1$ we have $c_1(1)=x_1$ directly from (\ref{GE0}). For any $k\ge 2$ by (\ref{Vkk}) we readily have
\be
v^{k,k}_{1,2_{k-1}}=(2t)^{k-1}.\label{vkk12}
\ee
By (\ref{ex2-gf}), (\ref{vkk12}), and (\ref{Vk2}), (\ref{initM}) we have: 
$$
c_2(1)=(v_{1,2}^{3,2}+v_{1,2}^{2,2})\big|_{t=0}x_1=(\partial_t (2t) +2t)\big |_{t=0}x_1=(2+2t)\big|_{t=0}x_1=2x_1
$$
thus the coefficient of $t^2$ is $2x_1/2=x_1$. 
By (\ref{Vks}):
$$
v^{4,3}_{1,2,2}=v^{3,2}_{1,2}(2t)+\partial_t v^{3,3}_{1,2,2}=\partial_t(2t) (2t)+\partial_t (2t)^2=2(2t)+4(2t) =12t
$$
therefore
$$
c_3(1)=(v_{1,2}^{4,2}+v_{1,2,2}^{4,3}+v_{1,2,2,2}^{4,4})\big|_{t=0}x_1=(\partial^2_t (2t) +12t+ (2t)^3)\big |_{t=0}x_1=0.
$$
which is the coefficient of $t^3$ in (\ref{TE}). Let us calculate just another coefficient. 
$$
c_4(1)=(v_{1,2}^{5,2}+v_{1,2,2}^{5,3}+v_{1,2,2,2}^{5,4}+v_{1,2,2,2,2}^{5,5})\big|_{t=0}x_1
$$
for which we need to calculate $v^{5,3}$, $v^{5,4}$:
\bea
&&v_{1,2,2,2}^{5,4}=v_{1,2,2}^{4,3}(2t)+\partial_t v_{1,2,2,2}^{4,4}=12t(2t)+\partial_t(2t)^3=24t^2+6(2t)^2=48t^2\nonumber\\
&&v_{1,2,2}^{5,3}=v_{1,2}^{4,2}(2t)+\partial_t v_{1,2,2}^{4,3}=\partial^2_t (2t)(2t)+\partial_t 12t=12\nonumber
\eea
and finally
$$
c_4(1)=(\partial_t^3(2t)+48t^2+12+(2t)^4)\big|_{t=0}x_1=12
$$
thus the coefficient of $t^4$ is $12x_1/ 4! = x_1/2$.\\

\subsubsection{The basic stationary Driver-type ODE} 
\label{sss:blons} 
$\dot x=ax^2$, for $a\in \re^+$. We have $\cal S=\{1\}$, $V=a$. The solution $(x,0)$ can be easily calculated by variable separation, 
the result is
\be
x(t)={x\over 1-axt}.\label{ex3-sol}
\ee
Let us calculate the solution through formula (\ref{cstat}): all coefficients  
$c_k(1)$ have each only one term in formula (\ref{Vkstat}) : $c_k(1)=v^{k+1}_{1,1_k}kv_{1,1}$, therefore, by successive substitutions:
$$
c_k(1)=v^{k+1}_{1,1_k}ka=v^k_{1,1_{k-1}}k(k-1)a^2=\ldots =k!a^k
$$
replacing which in (\ref{cstat}) yields
$$
x(t)=\sum_{k=0}^\infty c_k(1){t^k\over k!}=x\sum_{k=0}^\infty a^kx^kt^k={x\over 1-axt}
$$
which is (\ref{ex3-sol}). The convergence ray, $r$, is given by the ratio's criterion:
$$
{1\over r}=\lim_{k\to +\infty} {|a^{k+1}x^{k+1}|\over |a^kx^k|}=a|x|.
$$
As expected, the convergence interval is $(-1/a|x|, 1/a|x|)$: if $x>0$ (resp: $x<0$) the solution goes to $+\infty$ on the right  (resp: to $-\infty$ on the left)  end point. The limit function is defined and analytic for $t\le -1/a|x|$ (resp: for $t\ge 1/a|x|$), thus, by analytic prolongation, it is the maximal solution, defined on ${\bf I}_{x,0}=(-\infty, 1/a|x|)$ (resp; on ${\bf I}_{x,0}=(-1/a|x|,+\infty)$). \\

\subsubsection{The Airy equation}  
\label{sss:blons} 
$x^{(2)}=tx$. This is a second order linear ODE for which is not known a 'closed' form of the solution $(x,0)$ -- i.e. given through a finite composition of ordinary transcendental and/or elementary functions -- but the solution can be calculated by power series. The result\footnote{A way for calculating the solution, which can be found in many exercise books of analysis courses, 
is to set it equal to $\sum_{k=0}^\infty a_kt^k$ (the general form of a power series), and then calculating the second derivative (term to term). By substitution, equating the terms with the same degree, we obtain a recursive equation for the coefficients 
$a_{k+1}=a_{k-1} (k+1)^{-1}(k+2)^{-1}$, relating every coefficient to the coefficient placed two times before, which allows the calculation of all coefficients from $a_0$ and $a_1$. By setting $p_1=a_0$, $p_2=a_1$ we obtain (\ref{ex4-solf})} is a function 
of two parameters, namely $p_1,p_2$:
\be
x(t)=p_1+p_2t+{1\over 3!}p_1t^3+{2\over 4!}p_2t^4+{1\cdot 4\over 6!}p_1t^6+{2\cdot 5\over 7!}p_2t^7+{1\cdot 4\cdot 7\over 9!}p_1t^9+\ldots \label{ex4-solf}
\ee
and an application of the ratio criterion readily yields an infinite convergence ray. 
In order to find (\ref{ex4-solf}) by Theorem \ref{th:serie}, we need first to quadratize the Airy equation. Let us introduce the two indeterminates $x_1,x_2$ and set $x_1=\dot x$, $x_2=x$, then the solution $(x,0)$ of the Airy equation,  is the sub-solution $x_2$, of 
the solution $((x_1,x_2),0)$ of the linear system:
\be
\begin{array}{l}
\dot x_1=tx_2\\
\dot x_2=x_1
\end{array}
\label{ex4-tra}
\ee
Note that $(0,0)$ is an equilibrium point of (\ref{ex4-tra}), and the zero solution $x_1\equiv x_2\equiv 0$ is the only one singular solution. 
Therefore, all solutions, but the trivial one, are regular, thus (cf: \S \ref{ssc:sinreg}) they are equal to some solution $((x_1,x_2),0)$ with 
$x_1\ne 0$, and $x_2\ne 0$. 
Let us apply the Driver formula (\ref{dri}) to (\ref{ex4-tra}). We have two Driver variables 
$Z_{1,1}={x_2\over x_1}$ and $Z_{2,1}={x_1\over x_2}$. The equation for $Z_{1,1}$ is directly derived by the associated monomial 
$x_2=x_1^0x_2^1$ from which (using the exponents as coefficients): 
\be
\dot Z_{1,1}=(\pi^1_{1,1}v_1'Z_1 +\pi^1_{1,2} v_2'Z_2)Z_{1,1}=(-tZ_{1,1}+Z_{2,1})Z_{1,1}\label{eqZ11}
\ee
The equation for $Z_{2,1}$, in place of the Driver formula, can be directly derived observing that, by definition, $Z_{2,1}=Z^{-1}_{1,1}$, thus, it is the component $W_{1,1}$ of the inverse Driver (cf \S \ref{Inv-driv}):
\be
\dot Z_{2,1}=(tZ_{1,1}-Z_{2,1})Z_{2,1}.\label{eqZ21}
\ee
By renaming $x_3=Z_{1,1}$, $x_4=Z_{2,1}$, and gathering (\ref{ex4-tra}) -- rewritten as a final ODE --  (\ref{eqZ11}), and (\ref{eqZ21}) we 
obtain the (inclusive) quadratization:
\be
\begin{array}{l}
\dot x_1=tx_3x_1\\
\dot x_2=x_4x_2\\
\dot x_3=x_3x_2-tx_3^2\\
\dot x_4=tx_4x_3-x_4^2
\end{array}
\label{ex4-tra2}
\ee
whose frame is
\be
V=\left(
\begin{array}{cccc}
0&0&t &0\\
0&0&0&1\\
0&1&-t&0\\
0&0&t&-1
\end{array}
\right)
\label{ex4-fr}
\ee
Thus $\cal S=\{2,3,4\}$.
Now, let us verify that $x_2(t)$ calculated through (\ref{TE}) agrees with (\ref{ex4-solf}).  
We have $c_0(2)=x_2$, which is the zero-coefficient of (\ref{ex4-solf}) for $p_1=x_2$. As for $c_1(2)$:
$$
c_1(2)=\sum_{i_1\in \{2,3,4\}} v^{2,2}_{2,i_1} x_2x_{i_1}=v^{2,2}_{2,2} x^2_2 +v^{2,2}_{2,3} x_2x_{3} +v^{2,2}_{2,4} x_2x_{4}
$$
thus, recalling that $x_3={x_2\over x_1}=x_4^{-1}$, and observing, in (\ref{ex4-fr}), that $v_{2,2}=v_{2,3}=0$, $v_{2,4}=1$, we have:
$$
c_1(2)=x_1
$$
which agree with (\ref{ex4-solf}) provided that $x_2=p_1$. Now, the coefficient of $t^2$ in 
(\ref{ex4-solf}) is zero, so we expect to obtain  $c_2(2)|_{t=0}=0$. As a matter of fact:
\be
\begin{array}{rl}
c_2(2)&=\sum_{i_1\in \{2,3,4\}} v^{3,2}_{2,i_1}x_2x_{i_1} +\sum_{i_1,i_2\in \{2,3,4\}} v^{3,3}_{2,i_1,i_2} x_2x_{i_1}x_{i_2}\\
&=v^{3,2}_{2,2}x_2^2+v^{3,2}_{2,3}x_2x_{3} +v^{3,2}_{2,4}x_2x_{4}\\
&\qquad+ v^{3,3}_{2,2,2} x_2^3 +v^{3,3}_{2,2,3} x_2^2x_{3}+v^{3,3}_{2,2,4} x_2^2x_{4}\\
&\qquad+ v^{3,3}_{2,3,2} x^2_2x_3 +v^{3,3}_{2,3,3} x_2x^2_{3}+v^{3,3}_{2,3,4} x_2x_3x_{4}\\
&\qquad+ v^{3,3}_{2,4,2} x_2^2x_4 +v^{3,3}_{2,4,3} x_2x_4 x_{3}+v^{3,3}_{2,4,4} x_2x^2_{4}
\end{array}\label{ex4-tm}
\ee
from (\ref{Vk2}): $v^{3,2}_{2,i_1}=\partial_t v_{2,i_1}=0$,  for $i_1=2,3,4$,  and from (\ref{Vkk}) we have
\be
\begin{array}{l}
v^{3,3}_{2,2,i_2}=v_{2,2} (v_{2,i_2}+v_{2,i_2})=0, \quad i_2=2,3,4\\
v^{3,3}_{2,3,i_2}=v_{2,3} (v_{2,i_2}+v_{3,i_2})=0 \quad i_2=2,3,4\\
v^{3,3}_{2,4,2}=v_{2,4} (v_{2,2}+v_{4,2})=0\\
v^{3,3}_{2,4,3}=v_{2,4} (v_{2,3}+v_{4,3})=t\\
v^{3,3}_{2,4,4}=v_{2,4} (v_{2,4}+v_{4,4})=1-1=0
\end{array}\label{ex4-tm2}
\ee
therefore
$
c_2(2)=v^{3,3}_{2,4,3} x_2x_{4}x_{3}=t x_2
$
and $c_2(2)|_{t=0}=0$ as expected. 
Clearly, the more $k$ increases the more difficult is carrying out the calculation by hand 
for $c_k(2)$. It's still meaningful though, and not too cumbersome, to perform the calculation of  one more coefficient, $c_k(2)$, as the parameters $p_1,p_2$ have been fixed, and 
we expect $c_3(2)|_{t=0}=x_2$. We have
\be
c_3(2)=\sum_{i_1\in \{2,3,4\}} v^{4,2}_{2,i_1}x_2x_{i_1} +\sum_{i_1,i_2\in \{2,3,4\}} v^{4,3}_{2,i_1,i_2} x_2x_{i_1}x_{i_2}
+\sum_{i_1,i_2,i_3\in \{2,3,4\}} v^{4,4}_{2,i_1,i_2,i_3} x_2x_{i_1}x_{i_2}x_{i_3}
\label{ex4-c3}
\ee
and by (\ref{Vkk}): $v^{4,4}_{2,i_1,i_2,i_3}=v^{3,3}_{2,i_1,i_2}(v_{2,i_3} +v_{i_1,i_3}+v_{i_2,i_3})$. By (\ref{ex4-tm2}) 
$v^{3,3}_{2,4,3}=t$ is the only non zero element in the vector $v^{3,3}_2$, moreover, $v^{4,2}_{2,i_1}=\partial^2_tv_{2,i_1}=0$, 
and since $v^{3,2}_2=0$, by (\ref{Vks}) we have $v^{4,3}_2=\partial_t v^{3,3}_{2}$, which is zero but $v^{4,3}_{2,4,3}=1$.
Therefore (\ref{ex4-c3}) turns into
$$
c_3(2)=v^{4,3}_{2,4,3} x_2x_{4}x_{3}
+t\sum_{i_3\in \{2,3,4\}} (v_{2,i_3} +v_{4,i_3}+v_{3,i_3}) x_2x_{4}x_{3}x_{i_3}\\
$$
thus $c_3(2)|_{t=0}=v^{4,3}_{2,4,3} x_2x_{4}x_{3}=x_2$, as expected. Worth to be stressed is the fact that, the method of coefficients calculation we have pointed out in the latest note, applies just {\it after} having fixed the initial time-point ($t_0$, which has been set to zero), 
whereas, through the new method here presented we calculate the general expression, as a function of $t$, of $c_k(2)$ first, and then set $t=t_0=0$. Therefore we can use the same formulas for calculating the Taylor expansion with any different initial point $t_0$.\\

\section{Connection with the 22nd Hilbert's Problem\\ and Further Developments} 
\label{sc:HPFD}

\subsection{A partial solution to a differential form of the 22rd Hilbert's problem} 

The exact quadratization-based solution method described in the present paper can be interpreted as a kind of {\it uniformization} 
of a sub-class of analytic differential equation and as such it represents a (partial) solution, in the real case, 
of a {\it differential version} of the  
22nd Hilbert's problem \cite{Hi12}, as below explained. In the original formulation, the 22rd Hilbert's problem consists, given a set of complex polynomial (non differential) equations 
\be
\co^n\ni F(x,v)=0,\label{Fvx}
\ee
 in $n$ complex indeterminates, $x\in \co^n$,  with $m$ complex {\it parameters} (say $v$ the vector in $\co^m$ collecting them) 
in finding an open dense set in $\cal V\subset \co^m$ and a surjective analytic function $f:\cal V\to \co^n$, such that $f(v)$ solves (\ref{Fvx}), i.e. $(f(v),v)$ describes the set of solutions of (\ref{Fvx}) through $v$ (which means that this set is 'uniformized'). 
In the most general formulation, $x$ and $v$ are both 'indeterminates':  it is the uniformization $(f(v),v)$ that actually decides that $v$ is a set of 'parameters'. In other words: the uniformization is a parametrization of the solutions set such that (by the surjectivity of $f$) 
for any $x$ there is always a set of parameters $v$ such that $x$ solves (\ref{Fvx}). 
We do not consider here this more general formulation, since throughout this paper we have adopted, the more general {\it differential} setting on the one hand, but   
a more restrictive one on the other hand, in which all quantities are real, but 
we always know in advance what is a parameter and what an indeterminate, and that an unique solution always exists whatever we choose a set of parameters. 
A complete solution of the (non differential) uniformization problem has so far not been obtained, except for the one dimensional case, whereas the solution proposed 
by Poincair\'e in 1907 \cite{Po07} is considered a partial one. 

Now, provided we move to the real analytic case, and generalize the equations by allowing them 
1) to be ODEs:
\be
\re^n\ni G(w,x,\dot x,t)=0,\label{Fvxt}
\ee
depending on a set $w$ of real parameters (which can be even analytic functions of $t$) 
2)  for any fixed $w$, $G$ to belong to  the class of analytic $\sigma\pi$-reducible functions (thus, properly including the polynomials), of rank $n$ with respect to the $\dot x$ indeterminates, 
then Theorems \ref{th:serie}, and \ref{th:seriest}, provide a solution 
 of this generalized uniformization problem.  For,  note first that for any given $w$, a solution $(x,t_0)$ of (\ref{Fvxt}),  which is a real analytic function defined on some open interval ${\bf I}_{x,t_0}\subset \re$, can be identified with the sequence of coefficients 
 of its Taylor expansion: $c=\{c_k\}_{k=0}^\infty$, where $c_k=(c_k(1)(t_0),\ldots, c_k(n)(t_0))$. Let ${\cal S}_{t_0}$ the set 
 of all sequences $c$ such that the corresponding Taylor series centered at $t_0$ has a convergence ray $r>0$. Now, 
 as we have shown in \S \ref{sc:imp}, (\ref{Fvxt}) is quadratizable, thus $(x,t_0)$, for fixed $w$ is a sub-solution of some Driver in $d$ indeterminates 
 ($d\ge n$), and such a Driver is characterized by a frame $V(w)\in \cal V$ -- 
 where $\cal V$ denotes the set of all $d\times d$ matrices -- which is a function (analytic) of the original coefficients $w$. 
 Theorem \ref{th:serie} defines an analytic map $\cal V\to \cal S_{t_0}$, giving the solution in the form of a sequence $c\in \cal V$. 
 Therefore the quadratizing map combined with 
 the recursive equations (\ref{Vkk})--(\ref{Vk2}) defines a map $f:w\mapsto c$, such that $c=f(w)$ is the solution $(x,t_0)$ of 
 (\ref{Fvxt}). The latter is a kind of  generalized (to the differential case) uniformization for the ODE (\ref{Fvxt}). 
 
Such {\it differential} uniformization, clearly does not imply that the problem is solved as well in the non differential case, but surely 
that, if an uniformization exists for the set of solutions of the real analytic equation (\ref{Fvx}), for $F$ in the class of 
differentially $\sigma\pi$-reducible functions (i.e. such that the associated ODE is $\sigma\pi$-reducible), then the same uniformization is the differential uniformization of the ODE (\ref{Fvxt}), for $G(v,x,0,t)=F(x,v)$, thus of the set of all {\it constant solutions} of the ODE, i.e. of its {\it equilibrium points}. By the theory developed in the previous section, 
the latter problem amounts to finding the $x$, if any, such that 
$c_k=0$ $\forall k>0$ in (\ref{Vkk})--(\ref{Vk2}), so that $c_0=x$ and $x(t)\equiv x$, where $x(\cdot)$ is the curve of the solution 
$(x,t_0)$ of the ODE.\\

\subsection{Control problems, and differential flatness}  

We sketch here few guidelines for a possible development of the present research, by showing that the results of \S \ref{sc:SDTheo} 
can be used for finding a {\it flat output} for a control system. {\it Differential flatness} 
for control systems has been 
first introduced in \cite{FLMR92}, \cite{FLMR93}, in the context of differential algebra, and more recently defined in a more geometric 
Lie-B\"acklund framework relying on the paper \cite{Po93}. Most of the latter approach can be found in the textbook \cite{Le09}. We point out 
as well a further approach in which flatness can be described in terms of the notion of absolute equivalence defined by E. Cartan in \cite{Ca14}, 
\cite{Ca15b}, and where a nonlinear control system is viewed as a Pfaffian system on an appropriate space \cite{NiRM98}. 
The discussion below, mirrors the Lie-B\"acklund framework of \cite{Le09}, but it is translated in the setting of the present paper.

Let $y=y^{(0)}$, and $y^{(k)}$, $k$ integer, be a countable set of indeterminates in $\re^n$, and consider 
the following {\it countably infinite} ODEs:
\be
\partial_t y^{(k)}=y^{(k+1)}, \quad k\in \na.\label{inf-lin}
\ee
We denote by $({\bf y},t_0)$, where $\bf y$ is the sequence $\{y^{(k)}\}_{k\in \na}$, 
the generic solution of (\ref{inf-lin}) passing through $t_0\in \cal T\subset \re$, $\cal T$ an open interval, 
i.e. the function $t\mapsto {\bf y}(t)$, 
such that $y^{(k)}(t_0)=y^{(k)}$. Note that {\it any} map $t\mapsto {\bf y}(t)$, which is analytic on $\cal T$ in all components, 
determines a solution of (\ref{inf-lin}), where $y^{(k)}$ is the $k$th derivative of $y$ for any k. 
By reason of that (\ref{inf-lin}) is in fact said {\it the trivial system}. 
Now suppose that the ODE 
\be
\dot x=f(x,t),\label{Ofl}
\ee
is a Driver, with domain $\cal D\times \cal T\subset \re^{n+1}$, of some pre-Driver (i.e. it is a quadratization of some 
$\sigma\pi$-system, cf.   \S \ref{SD:sec}). In the most general situation, in turn, the latter $\sigma\pi$-system 
includes as sub-solutions all solutions of some original $\sigma\pi$-reducible system. 
Let $v_{i,j}$ the $n^2$ coefficients of (\ref{Ofl}), all defined and analytic on $\cal T$, and 
suppose there exist analytic functions $h_{i,j}: \cal U{\rm (open)}\subset \re^m \to \re^{n^2}$, with $m\le n^2$, such that 
$v_{i,j}=h_{i,j}(u)$, where 
$u\in \cal U$, is a vector of {\it controls}, analytic time-functions defined on $\cal T$, that can be arbitrarily chosen. 

A control problem
can be set as follows, given a function $t\mapsto x(t)$, find $u(t)$ such that $x(t)$ is the solution of (\ref{Ofl}) such that 
$x(t_0)=x$, for some $x\in \cal D$ and $t_0\in \cal T$. \footnote{We use here some abuse of terminology, in that what we call a  'control problem' is actually a kind   
of {\it inversion problem}, whereas what is actually of interest in control theory is the {\it tracking} of a trajectory $t\to x(t)$. 
Moreover,  we are considering just an {\it open loop} perspective.} 
Now, let us rewrite (\ref{Ofl}) by defining $g(x,u)=f(x,t)$, which makes explicit the dependence on $u$, define the new indeterminates $u^{(k)}$, $k\in \na$, with $u^{(0)}=u$, and consider the infinite dimensional ODE:
\be
\begin{array}{l}
\dot x=g(x,u),\\
\partial_t u=u^{(1)},\\
\partial_t u^{(1)}=u^{(2)}\\
\qquad \vdots\\
\end{array}\label{ext-sy}
\ee
whose solutions $((x,{\bf u}), t_0)$ -- $\bf u$ denoting the sequence $\{u^{(k)}\}_{k\in \na}$ -- include the solution $(x,t_0)$ of  (\ref{Ofl}) 
for $v_{i,j}=h_{i,j}(u)$ (and in fact it is the sub-solution of (\ref{ext-sy}) in the first $n$ components). 
We denote $\cal D^*$ (resp: $\cal Y$) the domain of (\ref{ext-sy}) (resp: of (\ref{inf-lin})) given by\footnote{$\re^n_{\infty}$ denotes the infinite product $\re^n\times \re^n \times \ldots$}:
\bea
&&\cal D^*=\cal D\times \cal U\times \re^m_{\infty},\label{domDst}\\
&&\cal Y=\re^n_\infty.\label{domcalY}
\eea
Note that, by the recursive equations (\ref{Vkk})--(\ref{Vk2}), there exist analytic functions $\phi^{k,s}$ such that 
\be
v^{k,s}=\phi^{k,s}(u,u^{(1)},\ldots,u^{(k)}), \quad s=2,\ldots, k; \quad k\ge 2\label{vksu}
\ee
then Lemma \ref{lm:serie} and Theorem \ref{th:serie} prove that there exists a map 
\be
{\Psi}: ((x,{\bf u}), t_0)\mapsto ({\bf y}, t_0),\label{mapPsi}
\ee 
such that $x=y$. As a matter of fact, such map is given by $y^{(k)}=q^{k+1}(x)$, $k\ge 0$, where $q^1(x)=x$, and for $k\ge 2$, 
the $q^k$'s are the polynomials 
given recursively by (\ref{qkrec}). Then, formula (\ref{qkrec}), the definition of $A_i$ given in (\ref{dA}), and the differential formula 
(\ref{papk}), shows that $p^{k+1}$ is a polynomial in $x$ and in $v, \partial_t v, \ldots \partial_t^k v$. Therefore, by (\ref{vksu}) for any 
$k\in \na$ there exist, {\it and can be calculated}, analytic functions $\psi^k$, such that  
\be
y^{(k)}=\psi^k(x,u,u^{(1)},\ldots,u^{(k)}) \label{ykpsi}
\ee
and the map $\Psi$ in (\ref{mapPsi}) is given by the infinite aggregate of the maps (\ref{ykpsi}).
The control problem is solved, as soon as one is able to show that the map (\ref{mapPsi}) is invertible, and, in particular, it is a  
diffeomorphism from $\cal D^*$ onto its own image in $\cal Y$. Here we are assuming a setting similar as in \cite{Le09}, and  we consider 
$\cal D^*$ and $\cal Y$ as differentiable manifolds (in infinite dimension), with the underlying topology generated by sets of the type: 
$U\times \re^m_\infty$ -- for the manifold (\ref{domDst}) --  with $U$ open set in $\cal D\times (\re^m)^N$ (for some positive integer $N$),  
and $U\times \re^n_\infty$ -- for the manifold (\ref{domcalY}) -- with $U$ open set in $\cal D\times (\re^n)^N$. 
A diffeomorphism $\Phi$ between manifolds endowed with the topology defined above, is said a {Lie-B\"acklund isomorphism}: it has components depending on a finite set of arguments only \cite{Le09}, thus, there exist analytic maps $\alpha,\beta$, and an 
integer $\rho$, such that
\bea
&&x=\alpha(y,y^{(1)}, \ldots, y^{(\rho)})\label{x-alp}\\
&&u=\beta(y,y^{(1)}, \ldots, y^{(\rho)})\label{x-bet}
\eea
Note that (\ref{x-alp}) in our setting is always verified by $x=y$, which is given by the inverse of 
the first $n$ components of the map $\Psi$ in (\ref{mapPsi}) (which is the identity in $\re^n$), thus, the solution to the control problem 
in open loop, amounts to  finding $\beta$ in (\ref{x-bet}) only.  Generally speaking, Theorem \ref{th:serie} suggests a principle way 
for calculating $\beta$: 
noticing that $c_k(i)=y_i^{(k)}$, one has to, first,  find an integer $\rho$, and $\rho$ equations among (\ref{GE0}), (\ref{GEs}), 
that are solvable with respect to a finite number, upper bounded by $\rho$, of
$v^{k,s}$, and next 
in solving -- on account of (\ref{vksu}) --  a finite set of differential equations, extracted from (\ref{Vkk})--(\ref{Vk2}), with respect to $u$. 
It should be noted that, a necessary condition for the map $\Phi$ to be a Lie-B\"acklund isomorphism, is that the flat output $y$ has the same dimension of the input $u$ \cite{Le09}. This is verified in our scheme, where $x=y$, only in a situation (not usual for control systems) 
where $m=n$. As a matter of fact, since (\ref{Ofl}) is a self-Driver, the {\it original} system has a lower number of indeterminates, say $l\le n$, 
and since (\ref{Ofl}) is a quadratization of the original system, there are just $l$ out of the $n$ $x_i$'s which we are interested to, and typically 
$m\le l$ (the number of controls is lower than the number of states). Nevertheless, Lemma \ref{lm:serie} holds {\it component-wise}, therefore 
we can even set $y\in \re^m$, and $y^{(k)}_{i_s}=q_{i_s}^{k+1}$, and obtain a trivial system as well, for a suitable choice of $i_1,\ldots,i_m$ such that $y^{(k)}_{i_s}=x^{(k)}_{i_s}$, $\forall k\in \na$, and $x_{i_s}$, $s=1,\ldots, m$ are $m$ out of the $l$ indeterminates associated to the solution of the original system. 
By using this setting, the jet-manifold (\ref{domcalY}), i.e. the domain of the trivial-system, 
becomes $\cal Y=\re^m_{\infty}$, 
the necessary condition for a Lie-B\"acklund isomorphism can be verified, $\alpha$ in (\ref{x-alp}) is no more an identity, 
but the problem is well-posed in any realistic situation. 

We highlight that, as noted in \cite{Le09}, flatness can be seen as a solution of another kind of generalized (in the framework of manifolds of jets of infinite order) 
uniformization problem of the 22rd Hilbert's problem.

\section{Conclusion} 

The goal of this paper was twofold: first, to give new results (generalizations and/or improvements) and further insight on the topic of exact quadratization, first issued in the article \cite{Ca15}, and second, which is the main result, to show that exact quadratization  allows to write the solution of any analytic quadratizable ODE (roughly speaking: 
{\it practically all} analytic ODEs) as a power series, directly from the ODE coefficients. 

Sections I and II has been devoted to the first task. In \S I we have summarized the basic result of \cite{Ca15}, but using a quite different  terminology, which in our intention aimed to focus the {\it strictly essential} features of {\it quadratizable} ODEs. In \S   \ref{ssc:sinreg}
we have given the definitions, of singular and regular part of an analytic system. which are used to prove one of the two most important new results of the section, 
which is Theorem \ref{th:GEQ} showing the globality of the exact quadratization of $\sigma\pi$-systems. Theorem \ref{th:GEQ} relies on 
Theorem \ref{th:EQns}, which is a re-formulation,  
in the setting of the present paper -- and where also a shorter proof is offered --  of the basic result on exact quadratization 
for $\sigma\pi$-systems of \cite{Ca15}.  The second most relevant result of the section is given in Theorem \ref{th:imp}, which is the generalization of exact quadratization to implicitly defined $\sigma\pi$-ODEs. A lot of top-up issues on exact quadratization -- that were not given or not enough scrutinized in \cite{Ca15} --  are given throughout the rest of the section, such as: definitions and properties of  the {\it inverse Driver} (\S \ref{Inv-driv}) which allows to view  
exact quadratization as a multi-variable generalization of  the {\it Bernoulli-type} ODEs ({\it Ibid.}), as well as the concept of inclusive quadratization (Definition \ref{df:inc-qua}), Driver-type ODE, and self-Driver (given in \S \ref{SD:sec}), 

In \S III we have given the main result of the paper, stated in  Theorem \ref{th:serie}. 
The general recursive formula (\ref{Vkk})--(\ref{Vk2})  allows the calculation, directly from the coefficients $v_{i,j}$,  
of the coefficients $c_k(i)$, for any $k\in\na$, of the Taylor expansion (\ref{SPSD}) of the $i$th component $x_i(t)$ of the 
solution of any {\it Driver-type} ODE,  
and, hence, by Theorem \ref{th:Drty}, and by the exact quadratization theorem, of the solution of any $\sigma\pi$-reducible ODE. 
Theorem \ref{th:seriest} yields the stationary case, in which the recursive formula for the coefficients calculation, simplifies in 
(\ref{Vkstat}). Theorem  \ref{th:serie}.is proven on the basis of a preliminary result, given in Lemma \ref{lm:serie}, which 
has an importance in its own right, since formula (\ref{kder}) allows to write any $\sigma\pi$-reducible ODE as a {\it trivial system}, 
a fact that we have highlighted in the final section, \S  \ref{sc:HPFD},  where the interplay of exact quadratization 
with the 22rd Hilbert's problem has been discussed, 
as well as with the concept of {\it differential flatness}, used in control theory.

%\appendix
%
%\section{Proofs of Preliminary Results}
%

% \begin{figure}[h]
% \centering
% %\psfrag{xaxis}{$n_1$} \psfrag{yaxis}{$\bar v$}
% \includegraphics[width=\columnwidth]{treasure_nMC_250}\caption{Numerically
% determined values of the a-posteriori outcome $\bar v$ (cf.
% Section~\ref{sec:apost}) for different values of $n_{1}$. In these
% experiments, the number of points is $N=10$, side length of the square region
% is 50 units, $m_{1} = \bar n_{2} = 10$, $\delta= 0.01$, $\beta= 10^{-5}$, and
% the rows and the columns were drawn uniformly randomly. }%
% \label{fig:treasure}%
% \end{figure}

% \begin{figure}[h]
% \centering
% %\psfrag{xaxis}{$n_1$} \psfrag{yaxis}{$\bar v$}
% \includegraphics[width=\columnwidth]{treasure_n2_1000}\caption{Numerically
% determined values of the a-posteriori outcome $\bar{v}$ (cf.
% Section~\ref{sec:apost}) for different values of $n_{1}$. In these
% experiments, $N=10$, $m_{1}=10,\bar{n}_{2}=1000$, $\delta=0.01$,
% $\beta=10^{-5}$, and the rows and the columns were drawn uniformly randomly.}%
% \label{fig:treasure1000}%
% \end{figure}

\end{document}